\documentclass[10pt]{amsart}

\usepackage{amsmath}
\usepackage{amssymb}
\usepackage{bm}
\usepackage{graphicx}
\usepackage{psfrag}
\usepackage{color}
\usepackage{hyperref}
\hypersetup{colorlinks=true, linkcolor=blue, citecolor=magenta, urlcolor=wine}
\usepackage{url}
\usepackage{algpseudocode}
\usepackage{fancyhdr}
\usepackage{mathtools}
\usepackage{tikz-cd}
\usepackage{xy}
\input xy
\xyoption{all}
\usepackage{stmaryrd}
\usepackage{calrsfs}

\newtheorem*{maintheorem*}{Main Theorem}
\newtheorem{theorem}{Theorem}[section]
\newtheorem*{theorem*}{Main Theorem}

\newtheorem{question}[theorem]{Question}
\newtheorem{prop}[theorem]{Proposition}

\newtheorem{lemma}[theorem]{Lemma}
\newtheorem{cor}[theorem]{Corollary}

\theoremstyle{definition}
\newtheorem{definition}[theorem]{Definition}
\newtheorem{remark}[theorem]{Remark}
\newtheorem{example}[theorem]{Example}
\numberwithin{equation}{section}

\newcommand{\cc}{\mathbb{C}}
\newcommand{\ff}{\mathbb{F}}
\newcommand{\nn}{\mathbb{N}}
\newcommand{\pp}{\mathbb{P}}
\newcommand{\qq}{\mathbb{Q}}
\newcommand{\rr}{\mathbb{R}}
\newcommand{\zz}{\mathbb{Z}}

\newcommand{\uu}{\mathcal{U}}

\providecommand\ldb{\llbracket}
\providecommand\rdb{\rrbracket}

\newcommand{\gp}{\text{gp}}

\newcommand{\supp}{\text{supp}}
\newcommand{\ii}{\mathcal{A}}
\newcommand{\gd}{\text{gd}}



\begin{document}
	
\mbox{}
\title{Divisibility and a weak ascending chain condition on principal ideals}

\author{Felix Gotti}
\address{Department of Mathematics\\MIT\\Cambridge, MA 02139}
\email{fgotti@mit.edu}

\author{Bangzheng Li}
\address{Department of Mathematics\\MIT\\Cambridge, MA 02139}
\email{liben@mit.edu}

\date{\today}
	
\begin{abstract}
	An integral domain $R$ is atomic if each nonzero nonunit of $R$ factors into irreducibles. In addition, an integral domain $R$ satisfies the ascending chain condition on principal ideals (ACCP) if every increasing sequence of principal ideals (under inclusion) becomes constant from one point on. Although it is not hard to verify that every integral domain satisfying ACCP is atomic, examples of atomic domains that do not satisfy ACCP are notoriously hard to construct. The first of such examples was constructed by A. Grams back in 1974. In this paper we delve into the class of atomic domains that do not satisfy ACCP. To better understand this class, we introduce the notion of weak-ACCP domains, which generalizes that of integral domains satisfying ACCP. Strongly atomic domains were introduced by D. D. Anderson, D. F. Anderson, and M. Zafrullah in 1990. It turns out that every weak-ACCP domain is strongly atomic, and so we introduce a taxonomic classification on our class of interest: ACCP implies weak-ACCP, which implies strong atomicity, which implies  atomicity. We study this chain of implications, putting special emphasis on the weak-ACCP property. This allows us to provide new examples of atomic domains that do not satisfy ACCP.
\end{abstract}

\bigskip
\maketitle

\section{Introduction}
\label{sec:intro}

Let $R$ be an integral domain. Following P. Cohn~\cite{pC68}, we say that $R$ is atomic if each nonzero nonunit of $R$ factors into irreducibles. Various relevant classes of integral domains are atomic, including Noetherian domains and Krull domains. Although the notion of atomicity in integral domains has shown up in the literature since the sixties, the first systematic study of atomic domains was carried out by D. D. Anderson, D. F. Anderson, and M. Zafrullah in~\cite{AAZ90}, where a taxonomic diagram to study factorizations in atomic domains was introduced. In the same paper, one of the notions considered by the authors is that of strongly atomic domains. Following their terminology, we say that $R$ is strongly atomic if any $r,s \in R$ have a common divisor $d$ in $R$ such that~$d$ factors into irreducibles and $\gcd(\frac rd, \frac sd) = 1$. Observe that 
every strongly atomic domain is atomic.
\smallskip

The integral domain $R$ is said to satisfy the ascending chain condition on principal ideals (ACCP) if every increasing sequence of principal ideals (under inclusion) becomes constant from one point on. Back in the sixties, it was wrongly asserted by Cohn in~\cite{pC68} that the property of being atomic and that of satisfying ACCP were equivalent. Although every integral domain satisfying ACCP is, in fact, strongly atomic \cite[Theorem~1.3]{AAZ90} and so atomic, not every (strongly) atomic domain satisfies ACCP. However, all known constructions of examples witnessing this statement are somehow sophisticated. The first complete construction was provided by A. Grams~\cite[Theorem~1.3]{aG74} back in 1974, although it is worth noticing that in 1973 Cohn himself suggested another construction~\cite[page~4]{pC73}\footnote{The example suggested by Cohn in~\cite{pC73} was kindly pointed out by S. Tringali.}. Grams' construction has been recently generalized by the authors in~\cite{GL23}. Further constructions of atomic domains that do not satisfy ACCP have been given by A. Zaks~\cite[Theorem~1.1 and Example~2]{aZ82} and M. Roitman~\cite{mR93}, and more recently, by J. Boynton and J. Coykendall~\cite[Example~3.5]{BC19} and also by the authors~\cite[Proposition~3.6]{GL22}. The connection of atomicity and the ACCP property has been recently explored in~\cite{BBNS23} in the non-commutative setting, where J. Bell, K. Brown, Z. Nazemian, and D. Smertnig have constructed a non-commutative finitely presented semigroup ring that is atomic but does not satisfy the ascending chain condition on principal right (or left) ideals.
\smallskip

In this paper we introduce the notion of weak-ACCP. We say that $R$ is weak-ACCP if it is atomic and every nonempty finite subset $S$ of $R \setminus \{0\}$ has a common divisor $d$ such that $\frac sd$ satisfies ACCP for some $s \in S$ (an element $x \in R$ satisfies ACCP if every ascending chain of principal ideals starting at $Rx$ becomes constant). Notice that 
every integral domain satisfying ACCP is a weak-ACCP domain. There are, however, weak-ACCP domains that do not satisfy ACCP (see Example~\ref{ex:weak-ACCP not ACCP} and Theorem~\ref{thm:weak-ACCP monoid ring not ACCP}). On the other hand, it turns out that every weak-ACCP domain is strongly atomic (see Proposition~\ref{prop:classes between atomicity and ACCP}). Finally, \cite[Example~5.2]{mR93} and Corollary~\ref{cor:weak-ACCP and polynomials}, used in tandem, yield an example of a strongly atomic domain that is not a weak-ACCP domain. Thus, based on some of the results and examples provided in Section~\ref{sec:weak-ACCP}, we obtain that the chain of implications [ACCP $\Rightarrow$ weak-ACCP $\Rightarrow$ strong atomicity $\Rightarrow$ atomicity] provides a proper taxonomic classification for the classes of atomic domains and atomic monoids.
\smallskip

In Section~\ref{sec:algebraic constructions}, we study the weak-ACCP property through the lenses of three algebraic constructions: the $D+M$ construction, localization, and directed unions. Both atomicity and the ACCP property were considered with respect to the $D+M$ construction in~\cite[Proposition~1.2]{AAZ90}. We offer a parallel result for the weak-ACCP property in Proposition~\ref{prop:D+M construction}. In addition, both atomicity and the ACCP property were studied with respect to localization in~\cite{AAZ92}, where the authors found conditions on a multiplicative subset $S$ of an integral domain $R$ such that these two properties (and further factorization properties) transfer between the integral domains $R$ and $R_S$ \cite[Theorems~2.1 and~3.1]{AAZ92}. In the direction of \cite{AAZ92}, we establish similar results for the weak-ACCP property in Propositions~4.5 and~4.6. Finally, we show that the weak-ACCP property transfers from the members of a directed family of integral domains to their directed union provided that every extension in the directed family is inert. 
\smallskip

In Section~\ref{sec:monoid rings}, we restrict our attention to the class of monoid rings. First, we prove that for any reduced torsion-free monoid $M$ satisfying ACCP, the weak-ACCP property transfers from any integral domain $R$ to the monoid ring $R[M]$ of $M$ over $R$. A similar result for the ACCP property has been given by D. D. Anderson and J. R. Juett in~\cite[Theorem~13]{AJ15} and previously observed by R. Gilmer and T. Parker in~\cite[Section~7]{GP74}. In Theorem~\ref{thm:Grams' construction and the weak-ACCP}, we provide a generalization of the celebrated Grams' construction of an atomic domain that does not satisfy ACCP. We include two applications of this theorem: one of them yields new examples of atomic domains that do not satisfy ACCP and the other one offers some insight of the behavior of atomicity under intersection. We conclude the section with a construction of a monoid ring that is weak-ACCP but does not satisfy ACCP (Theorem~\ref{thm:weak-ACCP monoid ring not ACCP}). This construction improves upon a construction given by the authors in \cite[Theorem~4.4]{GL23} as it lowers the rank of the monoid of exponents, which was infinite in the previous construction.
\smallskip

In~\cite{BC19}, Boynton and Coykendall exhibited an example of an atomic domain that does not satisfies ACCP. The arguments to establish that this example satisfies the desired properties make use some theoretical results about certain pullback construction that they developed in the first part of the same paper. In Section~\ref{sec:pullback}, we prove that the pullback domain provided by Boynton and Coykendall is indeed a weak-ACCP domain, refining the fact that this pullback is an atomic domain. It is worth emphasizing that we follow a completely different approach, which does not used the known fact the same pullback domain is atomic.

\bigskip
\section{Preliminary}
\label{sec:prelim}

\smallskip
\subsection{General Notation}

As it is customary, $\zz$, $\qq$, $\rr$, and $\cc$ will denote the set of integers, rational numbers, real numbers, and complex numbers, respectively. We let $\nn$ and $\nn_0$ denote the set of positive and nonnegative integers, respectively. In addition, we let $\pp$ denote the set of primes. For $p \in \pp$ and $n \in \nn$, we let $\ff_{p^n}$ be the finite field of cardinality $p^n$. For $a,b \in \zz$ with $a \le b$, we let $\ldb a,b \rdb$ denote the set of integers between $a$ and $b$, i.e., $\ldb a,b \rdb = \{n \in \zz \mid a \le n \le b\}$. In addition, for $S \subseteq \rr$ and $r \in \rr$, we set $S_{\ge r} = \{s \in S \mid s \ge r\}$ and $S_{> r} = \{s \in S \mid s > r\}$.

\smallskip
\subsection{Commutative Monoids}

A \emph{monoid} is a semigroup with an identity element. However, in the context of this paper we will tacitly assume that monoids are both cancellative and commutative. Let~$M$ be a monoid (additively written). We say that $M$ is \emph{torsion-free} provided that for all $b,c \in M$, if $nb = nc$ for some $n \in \nn$, then $b=c$. The \emph{Grothendieck group} $\gp(M)$ of $M$ is the unique abelian group $\gp(M)$ up to isomorphism satisfying that any abelian group containing a homomorphic image of $M$ will also contain a homomorphic image of $\gp(M)$. The \emph{rank} of $M$ is the rank of the $\zz$-module $\gp(M)$, that is, the dimension of the $\qq$-space $\qq \otimes_\zz \gp(M)$. The rank-one torsion-free monoids that are not groups are, up to isomorphism, the nontrivial additive submonoids of $\qq_{\ge 0}$, and they have been systematically studied recently (see~\cite{CGG21,GGT19} and references therein). The group of invertible elements of $M$ is denoted by $\uu(M)$. The quotient monoid $M/\uu(M)$, denoted here by $M_{\text{red}}$, is the \emph{reduced monoid} of $M$, and~$M$ is \emph{reduced} if $\uu(M)$ is a trivial group, in which case $M_{\text{red}}$ can be canonically identified with~$M$. For $b,c \in M$, we say that $b$ \emph{divides} $c$ in $M$ if $c = b + b'$ for some $b' \in M$; in this case, we often write $b \mid_M c$. The monoid $M$ is a \emph{valuation monoid} if for all $b,c \in M$ either $b \mid_M c$ or $c \mid_M b$.
\smallskip

An element $a \in M \! \setminus \! \uu(M)$ is an \emph{atom} (or an \emph{irreducible}) if whenever $a = u + v$ for some $u,v \in M$, then either $u \in \uu(M)$ or $v \in \uu(M)$. The set of atoms of $M$ is denoted by $\mathcal{A}(M)$. The monoid $M$ is \emph{antimatter} if $\mathcal{A}(M)$ is empty. On the other hand, $M$ is \emph{atomic} if every non-invertible element can be written as a sum of atoms. A characterization of atomicity has been recently given in \cite[Corollary~2.6]{sT23}. A subset $I$ of $M$ is an \emph{ideal} of~$M$ if $I + M = I$ or, equivalently, $I + M \subseteq I$. The ideal~$I$ is \emph{principal} if $I = b + M$ for some $b \in M$. The monoid $M$ satisfies the \emph{ascending chain condition on principal ideals} (\emph{ACCP}) if every ascending chain of principal ideals of $M$ eventually stabilizes. It is well known and not hard to verify that every monoid that satisfies ACCP is atomic (see \cite[Theorem~3.10]{sT22} and \cite[Theorem~3.4]{CT23} for more general versions of this result).
\smallskip

\subsection{Factorizations}

Observe that 
$M$ is atomic if and only if its reduced monoid $M_{\text{red}}$ is atomic. Let $\mathsf{Z}(M)$ denote the free (commutative) monoid on $\mathcal{A}(M_{\text{red}})$, and let $\pi \colon \mathsf{Z}(M) \to M_\text{red}$ be the unique monoid homomorphism fixing every element of $\mathcal{A}(M_{\text{red}})$. If $z := a_1 \cdots a_\ell \in \mathsf{Z}(M)$ for some $a_1, \dots, a_\ell \in \mathcal{A}(M_{\text{red}})$, then $\ell$ is the \emph{length} of $z$ and is denoted by $|z|$. For every $b \in M$, we set
\[
	\mathsf{Z}(b) = \mathsf{Z}_M(b) = \pi^{-1} (b + \uu(M)).
\]
If $|\mathsf{Z}(b)| = 1$ for every $b \in M$, then $M$ is a \emph{unique factorization monoid} (\emph{UFM}). On the other hand, if~$M$ is atomic and $|\mathsf{Z}(b)| < \infty$ for every $b \in M$, then $M$ is a \emph{finite factorization monoid} (\emph{FFM}). It follows directly from the definitions that every UFM is an FFM. Now, for every $b \in M$, we set
\[
	\mathsf{L}(b) = \mathsf{L}_M(b) = \{ |z| \mid z \in \mathsf{Z}(b) \}.
\]
If $M$ is atomic and $|\mathsf{L}(b)| < \infty$ for every $b \in M$, then $M$ is a \emph{bounded factorization monoid} (\emph{BFM}). Notice that 
if a monoid is an FFM, then it is a BFM. In addition, it is not hard to argue that every BFM must satisfy ACCP \cite[Corollary~1.4.4]{GH06}.
\smallskip

Let $R$ be an integral domain. We let $R^* := R \setminus \{0\}$ and $R^\times$ denote the multiplicative monoid and the group of units of $R$, respectively. The integral domain $R$ is \emph{antimatter} (resp., \emph{atomic}) provided that the multiplicative monoid $R^*$ is antimatter (resp., atomic). Similarly, the factorization properties defined in the previous paragraph can be adapted to integral domains via their multiplicative monoids (indeed, all such atomic properties were first investigated in the special setting of integral domains). The integral domain~$R$ is a \emph{unique} (resp., \emph{finite}, \emph{bounded}) \emph{factorization domain} if $R^*$ is a unique (resp., finite, bounded) factorization monoid, respectively. Observe that this new definition of a unique factorization domain coincides with the standard definition of a UFD. As we do for monoids, for integral domains we use the corresponding acronyms UFD, FFD, and BFD. We write $\mathsf{Z}(R) = \mathsf{Z}(R^\ast)$ and, for every $r \in R^\ast$, we write $\mathsf{Z}(r) = \mathsf{Z}_{R^\ast}(r)$ and $\mathsf{L}(r) = \mathsf{L}_{R^\ast}(r)$. We let $\ii(R)$ denote the set of irreducibles of the integral domains $R$. Observe that $\ii(R) = \mathcal{A}(R^*)$.
\smallskip

\subsection{Monoid Domains}

For the rest of this section, we assume that the monoid $M$ is torsion-free. Following the terminology in~\cite{rG84}, we let $R[M] := R[X;M]$ denote the monoid ring of $M$ over $R$, that is, the ring consisting of all polynomial expressions with exponents in $M$ and coefficients in $R$. When~$R$ is a field, $R[M]$ is often called a \emph{monoid algebra} and, in this case, the \emph{rank} of $R[M]$ is defined to be the rank of the monoid $M$. Since~$M$ is torsion-free, it follows from \cite[Theorem~8.1]{rG84} that $R[M]$ is an integral domain. In addition, it follows from \cite[Theorem~11.1]{rG84} that
\[
	R[M]^\times = \{rX^u \mid r \in R^\times \text{ and } u \in \uu(M)\}.
\]
As $M$ is torsion-free and cancellative, \cite[Corollary~3.4]{rG84} ensures the existence of a total order in~$M$ compatible with its operation. Therefore we can assume that $M$ is a totally ordered monoid, and we do so. Let $f := c_n X^{q_n} + \dots + c_1 X^{q_1}$ be a nonzero polynomial expression in $R[M]$ for some coefficients $c_1, \dots, c_n \in R^*$ and exponents $q_1, \dots, q_n \in M$ satisfying $q_n > \dots > q_1$. Then $\deg f := \deg_{R[M]} f := q_n$ and $\text{ord} \, f := \text{ord}_{R[M]} \, f := q_1$ are the \emph{degree} and the \emph{order} of $f$, respectively. In addition, we call $\text{supp} \, f := \text{supp}_{R[M]} f := \{q_1, \dots, q_n\}$ the \emph{support} of $f$. Observe that when $M = \nn_0$, we recover the standard notions of degree, order, and support for polynomial rings.

\bigskip
\section{A Weaker Version of the ACCP Property}
\label{sec:weak-ACCP}

It can be easily checked that every integral domain satisfying ACCP is atomic. The converse of this statement was asserted without proof by P. M. Cohn~\cite{pC68} back in 1968. A few years later, A. Grams disproved the assertion by finding, for every field $F$, a monoid domain $F[M]$ whose localization at the multiplicative subset $\{f \in F[M] \mid f(0) = 0\}$ is atomic but does not satisfy ACCP. The monoid in Grams' example is
\begin{equation} \label{Grams monoid}
	M = \Big\langle \frac{1}{2^n p_n} \ \Big{|} \ n \in \nn \Big\rangle,
\end{equation}
where $(p_n)_{n \in \nn}$ is the strictly increasing sequence whose underlying set is $\pp$.

\begin{definition}
	We call the monoid $M$ in~\eqref{Grams monoid} the \emph{Grams monoid}.
\end{definition} More recently, the authors generalized Grams' construction introducing the notion of a greatest-divisor submonoids. Let $M$ be a monoid, and let $N$ be a submonoid of $M$. For each $b \in M$, a \emph{greatest divisor} of $b$ \emph{in}~$N$ is an element $d \in N$ satisfying the following two properties:
\begin{itemize}
	\item $d \mid_M b$ and
	\smallskip
	
	\item if $d' \mid_M b$ for some $d' \in N$, then $d' \mid_M d$.
\end{itemize}
Note that 
any two greatest divisors in $N$ of the same element of $M$ must be associates, and so if $M$ is reduced every element of $M$ has at most one greatest divisor in $N$. We say that $N$ is a \emph{greatest-divisor} submonoid of $M$ provided that every element of $M$ has a greatest divisor in $N$. When $M$ is reduced and $N$ is a greatest-divisor submonoid of $M$, for each $b \in M$ we let $\text{gd}_N(b)$ denote the unique greatest divisor of $b$ in $N$.

\begin{example} \label{ex:Grams monoid has a greatest-divisor submonoid}
	Let $M$ be the Grams monoid. Observe that $N = \langle \frac{1}{2^n} \mid n \in \nn_0 \rangle$ is a submonoid of~$M$. It was proved in \cite[Lemma~1.1]{aG74} that each element $b \in M$ can be uniquely decomposed as
	\begin{equation} \label{eq:Grams equation}
		b = q(b) + \sum_{n \in \nn} c_n \frac{1}{2^n p_n}
	\end{equation}
	provided that $q(b) \in N$ and $c_n \in \ldb 0, p_n - 1 \rdb$ for every $n \in \nn$ (all but finitely many $c_n$'s are zero). Using the uniqueness of the decomposition~\eqref{eq:Grams equation}, one can easily argue that $N$ is a greatest-divisor submonoid of $M$.
\end{example}

Since the submonoid $\big\langle \frac 1{2^n} \mid n \in \nn \big\rangle$ of the Grams monoid is a valuation monoid, Grams construction is a special case of the following result.

\begin{theorem}\cite[Theorem~3.3]{GL23} \label{thm:generalized Grams example}
	Let $F$ be a field, and let $M$ be an atomic reduced torsion-free monoid. Also, let~$N$ be a submonoid of~$M$ satisfying the following conditions:
	\begin{itemize}
		\item[(1)] $N$ is a valuation greatest-divisor submonoid of $M$, and
		\smallskip
		
		\item[(2)] $\mathsf{L}_M(m - \emph{gd}_N(m))$ is finite for every $m \in M$.
	\end{itemize}
	Then $F[M]_S$ is atomic, where $S = \{f \in F[M] \mid f(0) \neq 0\}$.
\end{theorem}

Motivated by the hypotheses of Theorem~\ref{thm:generalized Grams example}, we proceed to introduce a weaker version of the ACCP property in the setting of both atomic monoids and domains; we call it \emph{weak-ACCP}. As we will see in next sections, this weak notion of the ACCP property is preserved under various constructions, and this will allow us to find new classes of atomic domains that do not satisfy ACCP. For a monoid~$M$, we say that an element $b \in M$ \emph{satisfies ACCP} if every ascending chain of principal ideals starting at $b + M$ eventually stabilizes. The next lemma is a special case of \cite[Lemma~2.3]{sT23}. For the sake of completeness, we include here a short proof.

\begin{lemma} \label{lem:ACCP elements are atomic}
	Let $M$ be a monoid. Then every element of $M$ that satisfies ACCP can be written as a sum of atoms.
\end{lemma}

\begin{proof}
	Take $b_0 \in M$ such that $b_0$ satisfies ACCP. Suppose, by way of contradiction, that $b_0$ cannot be written as a sum of atoms. Hence $b = b_1 + c_1$ for some $b_1, c_1 \in M \setminus \uu(M)$ such that one of them, say~$b_1$, is not a sum of atoms in $M$. Similarly, we can write $b_1 = b_2 + c_2$ for some $b_2, c_2 \in M \setminus \uu(M)$ such that~$b_2$ is not a sum of atoms in $M$. Continuing in this fashion, we obtain sequences $(b_n)_{n \in \nn_0}$ and $(c_n)_{n \in \nn}$ whose terms belong to $M \setminus \uu(M)$ such that $b_n = b_{n+1} + c_{n+1}$ for every $n \in \nn_0$. Therefore $(b_n + M)_{n \in \nn_0}$ is an ascending chain of principal ideals in $M$. In addition, the fact that $b_n - b_{n+1} = c_{n+1} \notin \uu(M)$ for any $n \in \nn$ implies that $(b_n + M)_{n \in \nn_0}$ does not stabilize, contradicting that $b_0$ satisfies ACCP.
\end{proof}

\begin{prop} \label{prop:weak-ACCP}
	For a monoid $M$, the following statements are equivalent.
	\begin{enumerate}
		\item[(a)] $M$ is atomic and every nonempty finite subset $S$ of $M$ has a common divisor $d \in M$ such that $s-d$ satisfies ACCP for some $s \in S$.
		\smallskip
		
		\item[(b)] Every nonempty finite subset $S$ of $M$ has a common divisor $d \in \langle \mathcal{A}(M) \rangle$ such that $s-d$ satisfies ACCP for some $s \in S$.
	\end{enumerate}
\end{prop}

\begin{proof}
	(a) $\Rightarrow$ (b): This is obvious.
	\smallskip
	
	(b) $\Rightarrow$ (a): Take $c \in M \setminus \uu(M)$. By hypothesis, we can write $c = a_1 + \cdots + a_k + b$ for some $a_1, \dots, a_k \in \mathcal{A}(M)$ and $b \in M$ satisfying ACCP. Since $b$ satisfies ACCP, $b$ can be written as a sum of atoms by virtue of Lemma~\ref{lem:ACCP elements are atomic}, and so $c$ can also be written as a sum of atoms. Thus, $M$ is atomic.
\end{proof}

It is worth noting that a monoid $M$ may not be atomic even if every nonempty finite subset $S$ of~$M$ has a common divisor $d \in M$ such that $s-d$ satisfies ACCP for some $s \in S$. For instance, $(\qq_{\ge 0},+)$ is not atomic but, as it is a valuation monoid, for every nonempty subset $S$ of $\qq_{\ge 0}$, the element $\min S$ can play the role of both $d$ and $s$ in the described condition.

\begin{definition} \label{def:weak-ACCP}
	We say that a monoid $M$ is \emph{weak-ACCP} if $M$ satisfies any of the equivalent conditions in Proposition~\ref{prop:weak-ACCP}. An integral domain is \emph{weak-ACCP} if its multiplicative monoid is weak-ACCP.
\end{definition}

Notice that 
the weak-ACCP property falls in between the property of being atomic and that of satisfying ACCP; that is, the chain of implications [ACCP $\Rightarrow$ weak-ACCP $\Rightarrow$ atomic] holds. Moreover, none of these implications is reversible, as the next examples will corroborate.

If $M$ is an atomic monoid having a greatest-divisor submonoid $N$ satisfying conditions~(1) and~(2) in Theorem~\ref{thm:generalized Grams example}, then $M$ must be weak-ACCP. Indeed, for every finite subset $S$ of $M$ we can take $d := \min \{\gd_N(s) \mid s \in S\}$. Since $N$ is a valuation monoid, $d$ is a common divisor of $S$ and, after taking $s \in S$ satisfying $\gd_N(s) = d$, we see that $\mathsf{L}_M(s-d) = \mathsf{L}_M(s - \gd_N(s))$ is bounded, which implies that $s-d$ satisfies ACCP. This allows us to give examples of weak-ACCP monoids that do not satisfy ACCP. For such examples, and further examples throughout this paper, we appeal to \emph{positive monoids}, that is, additive submonoids of $\rr_{\ge 0}$ (positive monoids have been well-studied recently in connection to atomicity; see \cite{BCG21,fG19} and references therein).

\begin{example} \label{ex:Grams monoid is weak-ACCP}
	Fix $p \in \pp$, and then let $(p_n)_{n \ge 1}$ be the increasing sequence whose underlying set is $\pp \setminus \{p\}$. Now consider the positive monoid $G_p = \big\langle \frac{1}{p^n p_n} \ \big{|} \ n \in \nn_0 \big\rangle$, which is a slightly modified version of the Grams monoid. The submonoid $\big\langle \frac 1{p^n} \ \big{|} \ n \in \nn \big\rangle$ of $G_p$ is a valuation monoid satisfying condition~(2) of Theorem~\ref{thm:generalized Grams example}. Thus, $G_p$ is a weak-ACCP monoid. However, it does not satisfy ACCP because, for instance, the ascending chain of principal ideals $\big( \frac 1{p^n} + G_p \big)_{n \ge 1}$ does not stabilize.
\end{example}

\begin{example} \label{ex:rational cyclic semirings is weak-ACCP but not ACCP}
	Take $q \in \qq \cap (0,1)$ such that $q \neq \frac 1n$ for any $n \in \nn$, and consider the positive monoid $S_q =  \langle q^n \mid n \in \nn_0 \rangle$. The monoid $S_q$ is atomic (see \cite[Theorem~6.2]{GG18}), and it has been argued in~\cite[Example~3.6]{GL23} that the valuation monoid $\nn_0$ is a divisor-closed submonoid of $S_q$ satisfying condition~(2) of Theorem~\ref{thm:generalized Grams example}. Hence $S_q$ is a weak-ACCP monoid. On the other hand, $S_q$ does satisfy ACCP. Indeed, the identity $\mathsf{d}(q) q^n = (\mathsf{d}(q) - \mathsf{n}(q))q^n + \mathsf{d}(q) q^{n+1}$ holds for every $n \in \nn_0$, and so the sequence of principal ideals $(\mathsf{d}(q)q^n + \nn_0[q])_{n \in \nn_0}$ is ascending, even though it does not stabilize.
\end{example}

On the other hand, there are atomic monoids that are not weak-ACCP, as the following example illustrates.

\begin{example} \label{ex:strongly atomic PM not weak-ACCP}
	Let $(p_n)_{n \in \nn}$ be the strictly increasing sequence whose underlying set is $\pp$, and set $q_n = p_n p_{n+2}$ for every $n \in \nn$. Now consider the positive monoid $M = \big\langle \frac {1}{q_n} \mid n \in \nn \big\rangle$. It is not hard to verify that $M$ is atomic with set of atoms $\mathcal{A}(M) = \big\{ \frac{1}{q_n} \mid n \in \nn \big\}$. Note that $\big\{\frac 12, \frac 13 \big\}$ is a subset of~$M$. Let us argue now that $0$ is the only common divisor of $\frac 12$ and~$\frac 13$ in $M$. Suppose, by way of contradiction, that this is not the case, and take $m \in \nn$ such that $\frac 1{q_m}$ is a common divisor of $\frac 12$ and~$\frac 13$. We split the rest of the argument in the following two cases.
	\smallskip
	
	\textsc{Case 1:} $m$ is even. Since $\frac{1}{q_m} \mid_M \frac 12$ and $m$ is even, we can take $k, \ell \in \nn$ and $c_1, \dots, c_k, d_1, \dots, d_\ell \in \nn_0$ with $d_\ell \neq 0$ such that the following equality holds:
	\begin{equation} \label{eq:atomic not weak-ACCP}
		\frac 12 = \sum_{i=1}^k c_i \frac{1}{p_{2i-1} p_{2i+1}} + \sum_{i=1}^\ell d_i \frac{1}{p_{2i} p_{2i+2}}.
	\end{equation}
	Assume, in addition, that $\ell$ is the minimum positive integer such that \eqref{eq:atomic not weak-ACCP} holds. If $\ell = 1$, then after clearing denominators in \eqref{eq:atomic not weak-ACCP}, we can see that $p_2 p_4$ divides $d_1$, and so the right-hand side of \eqref{eq:atomic not weak-ACCP} is at least $1$, a contradiction. Therefore we assume that $\ell \ge 2$. In this case, after clearing denominators in \eqref{eq:atomic not weak-ACCP}, we see that $p_{2\ell + 2}$ divides $d_\ell$. As a consequence, we can write $d_\ell = d p_{2\ell + 2}$ for some $d \in \nn$ and, accordingly, we can rewrite \eqref{eq:atomic not weak-ACCP} as
	\[
		\frac 12 = \sum_{i=1}^k c_i \frac{1}{p_{2i-1} p_{2i+1}} + \bigg( \frac{d_{\ell -1} + d p_{2\ell - 2}}{p_{2\ell - 2} p_{2\ell}} + \sum_{i=1}^{\ell - 2} d_i \frac{1}{p_{2i} p_{2i+2}} \bigg),
	\]
	which contradicts the minimality of $\ell$.
	\smallskip
	
	\textsc{Case 2:} $m$ is odd. In this case, we use the fact that $\frac{1}{q_m} \mid_M \frac 13$, and we obtain a contradiction by mimicking the same argument given in Case~1.
	\smallskip
	
	Thus, the only common divisor of $\frac 12$ and $\frac 13$ in $M$ is $0$. So in order to argue that $M$ is not weak-ACCP, it suffices to show that neither $\frac 12$ nor $\frac 13$ is an ACCP element. To do so, note that for every $n \in \nn$, the equality
	\[
		\frac{1}{p_n} = \frac{1}{p_{n+2}} + (p_{n+2} - p_n) \frac{1}{p_n p_{n+2}}
	\]
	implies that $\frac1{p_{n+2}} \mid_M \frac 1{p_n}$ for every $n \in \nn$. Hence $\big( \frac{1}{p_{2n-1}} + M \big)_{n \in \nn}$ and $\big( \frac{1}{p_{2n}} + M \big)_{n \in \nn}$ are ascending chains of principal ideals of $M$ starting at $\frac 12 + M$ and $\frac 13 + M$, respectively, that do not stabilize. As a consequence, $M$ is an atomic monoid that is not weak-ACCP.
\end{example}

As confirmed by the previous examples, the weak-ACCP property falls strictly in between the property of being atomic and that of satisfying ACCP. Following~\cite[page~5]{AAZ90}, we say that a monoid~$M$ is \emph{strongly atomic} provided that for all $b,c \in M$ there is a common divisor $d \in \langle \mathcal{A}(M) \rangle$ of $b$ and $c$ in~$M$ such that the only common divisors of $b-d$ and $c-d$ in $M$ are units. In addition, an integral domain is \emph{strongly atomic} if its multiplicative monoid is strongly atomic. As the weak-ACCP property, the strongly atomic property falls between the property of being atomic and satisfying ACCP. Moreover, we will prove in Proposition~\ref{prop:classes between atomicity and ACCP} that every implication in the following diagram of atomic classes holds and that none of them is reversible.

\begin{equation} \label{diag:AAZ's atomic chain for monoids}
	\begin{tikzcd}[cramped]
		\textbf{ ACCP } \  \ \arrow[r, Rightarrow] \arrow[red, r, Leftarrow, "/"{anchor=center,sloped}, shift left=1.7ex] & \ \textbf{ weak-ACCP } \  \ \arrow[r, Rightarrow] \arrow[red, r, Leftarrow, "/"{anchor=center,sloped}, shift left=1.7ex]  & \textbf{ strongly atomic } \  \ \arrow[r, Rightarrow] \arrow[red, r, Leftarrow, "/"{anchor=center,sloped}, shift left=1.7ex] & \textbf{ atomic }
	\end{tikzcd}
\end{equation}
\smallskip

In order to prove the next proposition it is convenient to introduction the notion of a weak-GCD monoid. The monoid $M$ is \emph{weak-GCD} if any $b,c \in M$ can be written as $b = d+b'$ and $c = d+c'$ for some $d, b', c' \in M$ such that $\gcd(b', c') = 0$
\begin{prop} \label{prop:classes between atomicity and ACCP}
	For a monoid $M$, the following statements hold.
	\begin{enumerate}
		\item If $M$ is strongly atomic, then $M$ is atomic.
		\smallskip
		
		\item If $M$ is weak-ACCP, then $M$ is strongly atomic.
		\smallskip
		
		\item None of the implication in Diagram~\eqref{diag:AAZ's atomic chain for monoids} is reversible.
	\end{enumerate}
\end{prop}

\begin{proof}
	(1) This was proved in \cite[Theorem~1.3]{AAZ90} for integral domains, and the proof for monoids follows in a completely similar fashion.
	\smallskip
	
	(2) Let $M$ be a weak-ACCP monoid. Suppose, for the sake of a contradiction, that $M$ is not strongly atomic. Take elements $x,y \in M$ such that $\gcd(x-d, y-d) \neq 0$ for any common divisor $d \in \langle \mathcal{A}(M) \rangle = M$. Since $M$ is weak-ACCP, there is a common divisor $d_0$ of $x$ and $y$ in $M$ such that either $x - d_0$ or $y - d_0$ satisfies ACCP. As $\gcd(x - d_0, y - d_0) \neq 0$, there is a common divisor $a_1 \in \mathcal{A}(M)$ of $x - d_0$ and $y - d_0$ in $M$, which satisfies that $\gcd(x - d_0 - a_1, y - d_0 - a_1) \neq 0$. Now if $a_1, \dots, a_n \in \mathcal{A}(M)$ and $d_n = \sum_{i=1}^n a_i$ is a common divisor of $x$ and $y$ in $M$ satisfying that $\gcd(x - d_n, y - d_n) \neq 0$, then taking a common divisor $a_{n+1} \in \mathcal{A}(M)$ of $x - d_n$ and $y - d_n$ in $M$, we obtain an element $d_{n+1} := d_n + a_{n+1}$ that is a common divisor of $x$ and $y$ in~$M$ and satisfies that $x - d_n + M \subsetneq x - d_{n+1} + M$ and $y - d_n + M \subsetneq y - d_{n+1} + M$. Hence we can assume the existence of a sequence $(d_n)_{n \in \nn}$ (with $d_0 = 0$) of common divisors of $x$ and $y$ in $M$ such that none of the ascending chain of principal ideals $(x - d_n + M)_{n \in \nn_0}$ and $(y - d_n + M)_{n \in \nn_0}$ stabilize. This contradicts that $M$ is weak-ACCP.
\smallskip

	(3) The multiplicative monoid of the integral domain $A$ constructed in \cite[Example~5.1]{mR93} is atomic but not weak-GCD. By part~(1), a monoid is strongly atomic if and only if it is atomic and weak-GCD. Thus, the multiplicative monoid $A^*$ is an atomic monoid that is not strongly atomic. In addition, we have seen weak-ACCP monoids that do not satisfy ACCP in Examples~\ref{ex:rational cyclic semirings is weak-ACCP but not ACCP} and~\ref{ex:strongly atomic PM not weak-ACCP}. As a result, we only need to construct a strongly atomic monoid that is not weak-ACCP. 
	
	Let $M$ be the monoid in Example~\ref{ex:strongly atomic PM not weak-ACCP}. We proceed to prove that $M$ is strongly atomic. Fix $x,y \in M$, and let us find a common divisor $d$ of $x$ and $y$ in $M$ such that $\gcd(x-d, y-d) = 0$. To do this, we can first find a positive integer $N$ such that $x,y \in M' := \big\langle \frac{1}{q_n} \mid n \in \ldb 1,N \rdb \big\rangle$. Furthermore, we can assume that $N \ge 10$.
	\smallskip
	
	\noindent \emph{Claim 1.} If $\frac {1}{q_N} \nmid_{M'} x$, then $\frac{1}{q_n} \mid_M x$ guarantees that either $n \le N$ or $n \not\equiv N \pmod{2}$.
	\smallskip
	
	\noindent \emph{Proof of Claim 1.} Assume, by way of contradiction, that $n = N+2n_0$ for some $n_0 \in \nn$, then we can take $L \in \nn_0,\ell \in \nn$ and $c_1,c_2,\cdots,c_N,d_1,d_2,\dots,d_L,e_1,e_2,\cdots,e_\ell \in \nn_0$ with $d_L \ne 0,e_\ell \ne 0$ such that the following equality holds:
	\begin{equation} \label{eq:strongly atomic but not weak-ACCP}
		x=\sum_{i=1}^N c_i \frac 1{q_i} + \sum_{i=1}^L d_i \frac 1{q_{N+2i-1}} + \sum_{i=1}^\ell e_i \frac 1{q_{N+2i}}.
	\end{equation}
	Assume, in addition, that $L+\ell$ is minimal among all choices of $L\in\nn_0$ and $\ell\in\nn$ such that \eqref{eq:strongly atomic but not weak-ACCP} holds. We claim that $L=0$. If $L>0$, then after clearing denominators in \eqref{eq:strongly atomic but not weak-ACCP} we can see that $p_{N+2L+1} \mid d_L$, and so $d_L = dp_{N+2L+1}$. Then using the equation $d_L q_{N+2L-1}^{-1} = d p_{N+2L-3}$, we can rewrite~\eqref{eq:strongly atomic but not weak-ACCP} as
	\begin{equation*}
		x = \sum_{i=1}^N c_i \frac 1{q_i} + \left( \sum_{i=1}^{L-2} d_i \frac 1{q_{N+2i-1}} + (d_{L-1}+dp_{N+2L-3}) \frac 1{q_{N+2L-3}} \right) + \sum_{i=1}^\ell e_i \frac 1{q_{N+2i}},
	\end{equation*}
	which contradicts the minimality of $L+\ell$. So we can rewrite \eqref{eq:strongly atomic but not weak-ACCP} as follows:
	\begin{equation}\label{eq:rewrite strongly atomic but not weak-ACCP}
		x = \sum_{i=1}^N c_i \frac 1{q_i} + \sum_{i=1}^\ell e_i \frac 1{q_{N+2i}}.
	\end{equation}
	Now we split the rest of the proof of Claim~1 into the following two cases.
	\smallskip
	
	\textsc{Case 1:} $\ell=1$. In this case, after clearing denominators we can see that $p_{N+4} \mid e_1$, and so $e_1 = e p_{N+4}$. Then from the equation $e_1 \frac 1{q_{N+2}} = e p_{N+4} \frac 1{q_{N+2}} = e p_N \frac 1{q_N}$, we deduce that $\frac 1{q_N} \mid_{M'} x$, which contradicts our assumption.
	\smallskip
	
	\textsc{Case 2:} $\ell>1$. In this case, we can clear denominators to see that $p_{N+2\ell+2} \mid e_\ell$, from which we deduce the equality $e_\ell = ep_{N+2\ell+2}$. Now using that
	\[
		e_\ell \frac 1{q_{N+2\ell}} = ep_{N+2\ell+2} \frac 1{q_{N+2\ell}} = e p_{N+2\ell-2} \frac 1{q_{N+2\ell-2}}
	\]
	we can rewrite~\eqref{eq:rewrite strongly atomic but not weak-ACCP} as
	\begin{equation*}
		x = \sum_{i=1}^N c_i \frac 1{q_i} + \left( \sum_{i=1}^{\ell-2} e_i \frac 1{q_{N+2i}} + (e_{\ell-1} + ep_{N+2\ell-2}) \frac 1{q_{N+2\ell-2}} \right),
	\end{equation*}
	which contradicts the minimality of $\ell$.
	\smallskip
	
	\noindent \emph{Claim 2.} If $\frac{1}{q_{N-1}} \nmid_{M'} x$, then $\frac 1{q_n} \mid_M x$ guarantees that either $n \le N$ or $n \equiv N \pmod{2}$.
	\smallskip
	
	\noindent \emph{Proof of Claim 2.} It follows by mimicking the proof of Claim 1.
	\smallskip
	
	We are now in a position to find a common divisor $d$. To do so, set
	\[
		t_1 = \max \bigg\{ t \in \nn_0 \ \bigg{|} \ \frac t{q_N} \mid_{M'} x \ \text{ and } \ \frac t{q_N} \mid_{M'} y \bigg\}
	\]
	and
	\[
		t_2 = \max \bigg\{ t \in \nn_0 \ \bigg{|} \ \frac t{q_{N-1}} \mid_{M'} x - \frac{t_1}{q_N} \ \text{ and } \ \frac t{q_{N-1}} \mid_{M'} y - \frac{t_1}{q_N} \bigg\}.
	\]
	Now set $d_0 = t_1 \frac 1{q_N} + t_2 \frac 1{q_{N-1}}$, and note that $d_0 \mid_{M'} x$ and $d_0\mid_{M'} y$. In addition, we see that $\frac 1{q_N} \nmid_{M'} x-d_0$ or $\frac 1{q_N} \nmid_{M'} y-d_0$ and, similarly, $\frac 1{q_{N-1}} \nmid_{M'} x-d_0$ or $\frac 1{q_{N-1}} \nmid_{M'} y-d_0$. Now suppose that $\frac 1{q_n} \mid_M x-d_0$ and $\frac 1{q_n} \mid_M y-d_0$. Then Claims~1 and~2 guarantee that $n \le N$.
	
	Now if $\gcd(x-d_0, y-d_0) \neq 0$, then one can subtract from $x - d_0$ and $y - d_0$ a nonzero common divisor, and such a divisor must be at least $\frac 1{q_N}$. We can repeat this process as many times as possible (it cannot last forever since $\qq$ is Archimedean). As a result, we will obtain a common divisor~$d$ of $x$ and $y$ in $M$ such that $\gcd(x-d, y-d) = 0$. Hence we can conclude that $M$ is strongly atomic.
\end{proof}

We conclude this section with a few words about hereditary atomicity, a property stronger than atomicity that was recently studied by J. Coykendall, R. Hasenauer, and the first author in~\cite{CGH23} in connection to the ACCP property. An integral domain $R$ is called \emph{hereditarily atomic} provided that every subring of $R$ is atomic. Although every hereditarily atomic domain is atomic by definition, there are integral domains satisfying ACCP that are not hereditarily atomic.

\begin{example} \label{ex:ACCP not HA}
	Consider the polynomial ring $R = \qq[X]$. Since $R$ is a UFD, it must satisfy ACCP. On the other hand, we claim that the subring $S := \zz + X\qq[X]$ of $R$ is not atomic. Observe that no element of the form $qX$ with $q \in \qq \setminus \{0\}$ is irreducible in $S$: for instance, we can write $qX = 2(\frac q2 X)$ and neither $2$ nor $\frac q2 X$ belong to $S^\times$. Hence the element $X$ does not factor into irreducibles in $S$, and so $S$ is not atomic.
\end{example}

However, if $R$ is an integral domain satisfying ACCP such that $R^\times$ is the trivial group, then it follows from \cite[Proposition 2.1]{aG74} that every subring of $R$ satisfies ACCP and, therefore, $R$ must be hereditarily atomic. However, the following question still remains open.

\begin{question}
	Does every hereditarily atomic domain satisfy ACCP?
\end{question}

\bigskip
\section{Algebraic Constructions}
\label{sec:algebraic constructions}

\smallskip
\subsection{The $D+M$ Construction}

We proceed to discuss the weak-ACCP property in the context of the $D+M$ construction. The $D+M$ construction is a rich source of examples and counterexamples in commutative ring theory. Let $T$ be an integral domain, and let $K$ and~$M$ be a subfield of $T$ and a nonzero maximal ideal of $T$, respectively, such that $T = K + M$. For a subdomain $D$ of $K$, set $R = D + M$. This construction was introduced by Gilmer \cite[Appendix II]{rG68} in the setting of valuation domains, and it was later investigated by J. Brewer and E. Rutter~\cite{BR76} and by D. Costa, J. L. Mott, and M. Zafrullah~\cite{CMZ78} in a more general context.

With notation as in the previous paragraph, it is well known that $R$ is atomic (resp., satisfies ACCP) if and only if $T$ is atomic (resp., satisfies ACCP) and $D$ is a field \cite[Proposition~1.2]{AAZ90}. Since the weak-ACCP seats between atomicity and the ACCP property, one may expect to have a similar result for weak-ACCP. In the next proposition we argue that this is indeed the case. First, let us state the following lemma, whose proof is straightforward.

\begin{lemma} \label{lem:associates weak-ACCP sets}
	Let $M$ be an atomic monoid, let $S := \{s_1, \dots, s_k\}$ be a nonempty finite subset of~$M$, and let $S' := \{u_1 s_1, \dots, u_k s_k\}$ for some $u_1, \dots, u_k \in \uu(M)$. For each $i \in \ldb 1,k \rdb$, an element $d \in M$ is a common divisor of $S$ in $M$ such that $\frac{s_i}d$ satisfies ACCP if and only if $d$ is a common divisor of~$S'$ in $M$ such that $\frac{u_i s_i}d$ satisfies ACCP.
\end{lemma}

We are in a position to show that being weak-ACCP behaves well under the $D+M$ construction.

\begin{prop} \label{prop:D+M construction}
	Let $T$ be an integral domain, and let $K$ and $M$ be a subfield of $T$ and a nonzero maximal ideal of $T$, respectively, such that $T = K + M$. For a subdomain $D$ of~$K$, set $R = D + M$. Then~$R$ is weak-ACCP if and only if $T$ is weak-ACCP and $D$ is a field.
\end{prop}

\begin{proof}
	For the direct implication, suppose that $R$ is weak-ACCP. In particular, $R$ is atomic, which implies via \cite[Proposition~1.2]{AAZ90} that $T$ is atomic and $D$ is a field. To argue that $T$ is also weak-ACCP, let $S$ be a nonempty finite subset of $T$. In order to find a common divisor $d$ of $S$ in $T$ such that $\frac sd$ satisfies ACCP in~$T$ for some $s \in S$, Lemma~\ref{lem:associates weak-ACCP sets} allows us to replace each element of $S$ by any of its associates and, therefore, we can assume that $S \subseteq M \cup (1 + M)$. In particular, $S \subseteq R$. Since~$R$ is weak-ACCP, there is a common divisor $d$ of $S$ in~$R$ such that $\frac sd$ satisfies ACCP in $R$ for some $s \in S$. In particular, $d$ is a common divisor of~$S$ in~$T$. Thus, proving that $T$ is weak-ACCP amounts to arguing that $\frac sd$ satisfies ACCP in~$T$.
	
	Since $D$ is a field, every principal ideal of~$R$ (resp., of $T$) has one of the forms $R m$ or $R (1+m)$ (resp., $T m$ or $T (1+m)$) for some $m \in M$. Furthermore, for all $m_1, m_2 \in M$, one can readily verify the following statements:
	\begin{enumerate}
		\item $R (1+m_1) \subseteq R(1 + m_2)$ if and only if $T (1+m_1) \subseteq T (1 + m_2)$,
		\smallskip
		
		\item $R m_1 \subseteq R (1 + m_2)$ if and only if $T m_1 \subseteq T (1 + m_2)$,
		\smallskip
		
		\item $R m_1 \subseteq R m_2$ implies that $T m_1 \subseteq T m_2$, and 
		\smallskip
		
		\item $T m_1 \subseteq T m_2$ implies that $R m_1 \subseteq R k m_2$ for some $k \in K^\times$.
	\end{enumerate}
	Now let $(T r_n)_{n \in \nn}$ be an ascending chain of principal ideals of $T$ starting at $T \frac sd$ and assume, without loss of generality, that $r_n \in M \cup (1+M)$ for every $n \in \nn$. Moreover, by virtue of the statement~(4) and the inclusion $K^\times \subseteq T^\times$, we can assume that for each $j \in \nn$ with $r_j, r_{j+1} \in M$ the inclusion $R r_j \subseteq R r_{j+1}$ holds. Therefore it follows from statements~(1) and~(2) that $(R r_n)_{n \in \nn}$ is an ascending chain of principal ideals of~$R$ starting at $R \frac sd$, and so it must eventually stabilize. As a result, it follows from statements~(1) and~(3) that $(T r_n)_{n \in \nn}$ eventually stabilizes. Hence $\frac sd$ satisfies ACCP in $T$, which completes the direct implication.
	\smallskip
	
	For the reverse implication, suppose that $T$ is weak-ACCP and~$D$ is a field. Since $T$ is atomic and~$D$ is a field, $R$ must be atomic by \cite[Proposition 1.2]{AAZ90}. Now suppose that $S$ is a nonempty finite subset of~$R$ and assume, as we did to prove the direct implication, that $S \subseteq M \cup (1 + M)$. Since~$T$ is weak-ACCP, there is a common divisor~$d$ of $S$ in $T$ such that $\frac sd$ satisfies ACCP in $T$ for some $s \in S$. Because $s \in M \cup (1+M)$ and $D$ is a field, we can replace $d$ by one of its associates in $T$ so that both~$d$ and $\frac sd$ belong to $R$. Finally, since $\frac sd$ satisfies ACCP in~$T$, it follows from (1)--(4) above that every ascending chain of principal ideals of $R$ starting at $R \frac sd$ must stabilize, and so $\frac sd$ also satisfies ACCP in $R$. Hence $R$ is weak-ACCP.
\end{proof}

\smallskip
\subsection{Localization}

Neither atomicity nor the ACCP property are preserved under localization. Indeed, there are integral domains satisfying ACCP and having an antimatter localization. The following example sheds some light upon this observation.

\begin{example} \label{ex:ACCP and atomicity are not preserved under localization}
	Consider the monoid algebra $\ff_p[M]$, where $\ff_p$ is the field of $p$ elements  and $M$ is the positive monoid $\{0\} \cup \qq_{\ge 1}$. Observe that the localization of $\ff_p[M]$ at its multiplicative subset $S := \{X^m \mid m \in M\}$ is the group algebra $\ff_p[\qq]$. It follows from \cite[Theorem~14.17]{rG84} that $\ff_p[\qq]$ does not satisfy ACCP. Moreover, since $\ff_p$ is a perfect field of characteristic $p$, every nonunit element in the group algebra $\ff_p[\qq]$ is the $p$-th power of a nonunit element, and so $\ff_p[\qq]$ is an antimatter domain; that is, $\ff_p[\qq]$ does not contain irreducibles.
\end{example}

However, under certain condition on the localization extension, atomicity and the ACCP property are both preserved under localization. Let $A \subseteq B$ be a ring extension. Following Cohn~\cite{pC68}, we call $A \subseteq B$ an \emph{inert extension} if $xy \in A$ for $x,y \in B^\ast$ implies that $ux, u^{-1}y \in A$ for some $u \in B^\times$. Let $A \subseteq B$ be an inert extension of integral domains. Then one can readily verify that $\ii(A) \subseteq B^\times \cup \ii(B)$. Therefore if $A \subseteq B$ is inert and $A^\times = B^\times \cap A$, then $\ii(A) = \ii(B) \cap A$.

\begin{example} \hfill
	\begin{enumerate}
		\item 	If $R$ is an integral domain, then the extension $R \subseteq R[X]$ is inert.
		\smallskip
		
		\item  Furthermore, under the usual notation of the $D+M$ construction, it is not hard to argue that both extensions $D \subseteq R$ and $R \subseteq T$ are inert.
		\smallskip
		
		\item Fix $n \in \nn$ with $n \ge 2$, and then consider the extension $R[X^n] \subseteq R[X]$. Observe that 
		$R[X^n]^\times = R^\times$. Notice, on the other hand, that although $X^n \in R[X^n]$, there is no $u \in R^\times$ such that  $uX \in R[X^n]$. As a result, $R[X^n] \subseteq R[X]$ is not an inert extension.
	\end{enumerate}
\end{example}

For any integral domain $R$ with a multiplicative subset $S$ such that $R \subseteq R_S$ is an inert extension, it has been established in~\cite[Theorem~2.1]{AAZ92} that $R_S$ is atomic (resp., satisfies ACCP) provided that $R$ is atomic (resp., satisfies ACCP). As the next proposition indicates, the same holds for the weak-ACCP property.

\begin{prop} \label{prop:when R weak-ACCP implies R_S weak-ACCP}
	Let $R$ be an integral domain, and let $S$ be a multiplicative subset of $R$ such that the extension $R \subseteq R_S$ is inert. If $R$ is weak-ACCP, then $R_S$ is weak-ACCP.
\end{prop}

\begin{proof}
	Suppose that $R$ is weak-ACCP. Assume, without loss of generality, that $S$ is saturated, so $R_S^\times = S$.  Since $R$ is atomic and the extension $R \subseteq R_S$ is inert, $R_S$ is also atomic. To show that every nonempty finite subset $W$ of $R_S$ has a common divisor $d \in R_S$ such that $\frac wd$ satisfies ACCP in $R_S$ for some $w \in W$, by virtue of Lemma~\ref{lem:associates weak-ACCP sets} it suffices to assume that $W \subseteq R$. Write $W = \{w_1, \dots, w_n\}$. Since $R$ is weak-ACCP, we can choose a common divisor $d$ of $W$ in $R$ such that $\frac wd$ satisfies ACCP in~$R$ for some $w \in \{w_1, \dots, w_n\}$. Note that $d \mid_{R_S} w_i$ for every $i \in \ldb 1,n \rdb$. So we are done once we prove that every ascending chain on principal ideals $\big(R_S r_n)_{n \in \nn}$ of $R_S$ starting at $R_S \frac{w}{d}$ (with $r_n \in R$ for every $n \in \nn$) eventually stabilizes. As $\frac{w}{d}$ satisfies ACCP in $R$, the ascending chain of principal ideals $(R r_n)_{n \in \nn}$ eventually stabilizes. Now observe that if for some $n \in \nn$ the inclusion $R_S r_n \subsetneq R_S r_{n+1}$ holds, and so $r_n = r_{n+1} \frac rs$ for some nonzero nonunit $\frac rs \in R_S$, then the fact that $R \subseteq R_S$ is an inert extension ensures the existence of $u \in S$ with $u r_{n+1} \in R$ and $(us)^{-1}r \in R$; thus, $\frac{r_n}{r_{n+1}} = u(us)^{-1}r \in R \setminus R^\times$, which implies that $R r_n \subsetneq R r_{n+1}$. This, together with the fact that $(R r_n)_{n \in \nn}$ stabilizes, implies that $(R_S r_n)_{n \in \nn}$ also stabilizes. Hence $R_S$ is a weak-ACCP domain.
\end{proof}

Our next result is, in certain sense, the converse of Proposition~\ref{prop:when R weak-ACCP implies R_S weak-ACCP}. First, we need the following definition. A saturated multiplicative subset $S$ of an integral domain $R$ is called \emph{splitting} if every $r \in R$ can be written as $r = as$ for some $a \in R$ and $s \in S$ such that $Ra \cap Rs' = R a s'$ for all $s' \in S$. It has been proved in \cite[Theorem~3.1]{AAZ92} that if $S$ is a splitting multiplicative subset of $R$ generated by primes, then $R$ is atomic (resp., satisfies ACCP) provided that $R_S$ is atomic (resp., satisfies ACCP). We proceed to show that a similar statement holds for the weak-ACCP property.

\begin{prop}
	Let $R$ be an integral domain, and let $S$ be a splitting multiplicative subset of $R$ generated by primes. If $R_S$ is weak-ACCP, then $R$ is weak-ACCP.
\end{prop}

\begin{proof}
	Since $R_S$ is atomic, so is $R$ by \cite[Theorem~3.1]{AAZ90}. Let $N$ be the subset of elements of $R$ that are not divisible by any prime in $S$. Take $r \in R$ and write $r = bs$ for some $b \in R$ and $s \in S$ such that $R b \cap R s' = R bs'$ for all $s' \in S$. Let $p$ be a prime in $S$. The fact that  $R$ is atomic ensures the existence of a maximum $m \in \nn_0$ such that $p^m \mid_R b$, and so $b \in Rb \cap R p^m = R bp^m$, which implies that $m=0$. Therefore each $r \in R$ can be written as $r = bs$ for some $b \in N$ and $s \in S$. In addition, it is not hard to verify that $R_S b \cap R = R b$ for all $b \in N$ (see \cite[proof of Theorem~3.1]{GP74}).
	
	Let $W$ be a nonempty finite subset of $R$. Since $R_S$ is weak-ACCP, there is a common divisor $a$ of~$W$ in $R_S$ such that $\frac wa$ satisfies ACCP in $R_S$ for some $w \in W$, and we can assume that $a \in R$. As~$S$ is a splitting multiplicative subset of $R$, we can further assume that $R a \cap R s = R as$ for all $s \in S$. Take $w \in W$ and write $w = a \frac{r'}{s'}$ for some $r' \in R$ and $s' \in S$. Then $r'a \in Ra \cap Rs' = Ras' $ and so $s' \mid_R r'$. Therefore~$a$ is a common divisor of $W$ in $R$. Finally, let us argue that every ascending chain of principal ideals $(R t_n)_{n \in \nn}$ starting at $R \frac wa$ must stabilize. For every $n \in \nn$, write $t_n = b_n s_n$ for some $b_n \in N$ and $s_n \in S$ such that $R b_n \cap R s = R b_n s$ for all $s \in S$. For each $n \in \nn$, the inclusion $R t_n \subseteq R t_{n+1}$ implies that $R b_n \subseteq R b_{n+1}$ and $R s_n \subseteq R s_{n+1}$. Thus, $(R b_n)_{n \in \nn}$ and $(R s_n)_{n \in \nn}$ are also ascending chains of principal ideals of $R$. As a result, $(R_S b_n)_{n \in \nn}$ is an ascending chain of principal ideals of $R_S$ starting at $R_S b_1 = R_S \frac ra$ and, therefore, stabilizes. As $R_S b_n \cap R = R b_n$ for every $n \in \nn$, the chain $(R b_n)_{n \in \nn}$ also stabilizes. Since the chain $(R s_n)_{n \in \nn}$ stabilizes, we conclude that $(R t_n)_{n \in \nn}$ must stabilize as well. Hence $R$ is a weak-ACCP.
\end{proof}

\smallskip
\subsection{Directed Unions}

A partially ordered set (poset) $\Gamma$, whose order relation is denoted by $\preceq$, is \emph{directed} provided that for all $\alpha, \beta \in \Gamma$, there exists $\theta \in \Gamma$ such that $\alpha \preceq \theta$ and $\beta \preceq \theta$. Let $(\Gamma, \preceq)$ be a nonempty directed poset. A family $(R_\gamma)_{\gamma \in \Gamma}$ of integral domains indexed by $\Gamma$ is a \emph{directed family} of integral domains if for all $\alpha, \beta \in \Gamma$ with $\alpha \preceq \beta$ the integral domain $R_\alpha$ is a subring of $R_\beta$. If $(R_\gamma)_{\gamma \in \Gamma}$ is a directed family of integral domains, then $\bigcup_{\gamma \in \Gamma} R_\gamma$ is an integral domain, which is the \emph{directed union} of $(R_\gamma)_{\gamma \in \Gamma}$. In general, the directed union of atomic integral domains (resp., integral domains satisfying ACCP) need not be atomic (resp., satisfy ACCP). The following example illustrates this observation.

\begin{example}
	For every $n \in \nn$, consider the monoid algebra $\ff_2[M_n]$, where $M_n := \frac{1}{2^n} \nn_0$. Observe that, 
	for every $n \in \nn$, the monoid algebra $\ff_2[M_n]$ is isomorphic to $\ff_2[X]$ and, therefore, is a UFD. In particular, $\ff_2[M_n]$ satisfies ACCP for every $n \in \nn$. Observe now that $(\ff_2[M_n])_{n \in \nn}$ is a directed family of integral domains with directed union $\ff_2[M]$, where $M := \big\{ \frac{m}{2^n} \mid m,n \in \nn_0 \big\}$. Since~$\ff_2$ is a perfect field of characteristic $2$ and every element of $M$ is $2$-divisible (that is, $\frac q2 \in M$ when $q \in M$), each nonunit of $\ff_2[M]$ is the square of a nonunit. As a result, the monoid algebra $\ff_2[M]$ is antimatter, whence it is not even atomic.
\end{example}

It is proved in \cite[Theorem~5.2]{AAZ92} that if $(R_\gamma)_{\gamma \in \Gamma}$ is a directed family of integral domains that are atomic (resp., satisfy ACCP) such that $R_\alpha \subseteq R_\beta$ is an inert extension for all $\alpha, \beta \in \Gamma$ with $\alpha \preceq \beta$, then the directed union of the family $(R_\gamma)_{\gamma \in \Gamma}$ is atomic (resp., satisfies ACCP). The atomicity of the previous statement was first observed by Zaks in \cite{aZ82}. We proceed to show that a similar result holds for the weak-ACCP property.

\begin{prop}
	For a nonempty directed poset $(\Gamma, \preceq)$, let $(R_\gamma)_{\gamma \in \Gamma}$ be a directed family of integral domains such that the extension $R_\alpha \subseteq R_\beta$ is inert for any $\alpha, \beta \in \Gamma$ with $\alpha \preceq \beta$.  If $R_\gamma$ is a weak-ACCP domain for each $\gamma \in \Gamma$, then $\bigcup_{\gamma \in \Gamma} R_\gamma$ is a weak-ACCP domain.
\end{prop}

\begin{proof}
	Suppose that $R_\gamma$ is a weak-ACCP domain for each $\gamma \in \Gamma$, and set $R = \bigcup_{\gamma \in \Gamma} R_\gamma$. Since $R_\gamma$ is atomic for all $\gamma \in \Gamma$, it follows from \cite[Theorem~5.2]{AAZ92} that $R$ is atomic. Let $S$ be a nonempty finite subset of $R$. As $S$ is finite, we can take $\theta \in \Gamma$ such that $S \subseteq R_\theta$. Now the fact that $R_\theta$ is a weak-ACCP domain guarantees the existence of a common divisor $d \in R_\theta$ of $S$ in $R_\theta$ such that $\frac sd$ satisfies ACCP in~$R_\theta$ for some $s \in S$. Suppose that $(Rx_n)_{n \in \nn_0}$ is an ascending chain of principal ideals of $R$ with $x_0 = \frac sd$. Take $y_0 = \frac sd$, and suppose we have found $y_1, \dots, y_n \in R_\theta$ such that $y_k R = x_k R$ and $y_{k-1} R_\theta \subseteq y_k R_\theta$ for every $k \in \ldb 1,n \rdb$. Since $x_{n+1} \mid_R y_n$, the fact that $R_\theta \subseteq R$ is inert allows us to take $u_{n+1} \in R^\times$ such that $y_{n+1} := u_{n+1}x_{n+1} \in R_\theta$ divides $y_n$ in $R_\theta$, in which case $y_{n+1}R = x_{n+1} R$ and $y_n R_\theta \subseteq y_{n+1} R_\theta$. Observe, in addition, that the last inclusion is strict when $x_n R \subseteq x_{n+1} R$ is strict. Hence from the fact that $\frac sd$ satisfies ACCP in $R_\theta$, we deduce that $\frac sd$ satisfies ACCP in~$R$. Thus,~$R$ is a weak-ACCP domain.
\end{proof}

\bigskip
\section{Monoid Rings}
\label{sec:monoid rings}

In this section, we consider the weak-ACCP property through the lens of monoid rings. We will use the weak-ACCP property to generalize Grams' construction (see \cite{aG74}) so that we can construct further examples of atomic domains that do not satisfy ACCP. In addition, we will produce an atomic monoid ring of finite rank (rank~$2$) that does not satisfy ACCP (all known atomic monoid rings that do not satisfy ACCP have infinite rank).

\smallskip
\subsection{Ascent of the Weak-ACCP Property}

The weak-ACCP property transfers from an integral domain to any of its monoid rings provided that the monoid of exponents satisfies ACCP.

\begin{prop} \label{prop:weak-ACCP from monoid to monoid rings}
	Let $R$ be an integral domain, and let $M$ be a reduced torsion-free monoid satisfying ACCP. Then $R$ is weak-ACCP if and only if $R[M]$ is weak-ACCP.
\end{prop}

\begin{proof}
	For the reverse implication, suppose that $R[M]$ is weak-ACCP, and let $S$ be a nonempty finite subset of $R^*$. As $R[M]$ is weak-ACCP, there exists a common divisor $d \in R[M]$ of $S$ in $R[M]$ such that $\frac sd$ satisfies ACCP in $R[M]$ for some $s \in S$. Now the fact that $R^*$ is a divisor-closed submonoid of $R[M]^*$ guarantees that $d$ is a common divisor of $S$ in $R$ and also that $\frac sd$ satisfies ACCP in~$R$. Thus,~$R$ is weak-ACCP.
	\smallskip
	
	To prove the direct implication, suppose that $R$ is weak-ACCP. Since $M$ is torsion-free, $R[M]$ is an integral domain. Now let $S$ be a nonempty finite subset of $R[M]^*$. Let $C$ be the set of nonzero coefficients of all the polynomials in $S$. Since $C$ is a nonempty finite subset of the weak-ACCP domain~$R$, there exists a common divisor $d \in R$ of $C$ that can be written as a product of irreducibles and a polynomial expression $f = \sum_{i=1}^n c_i X^{m_i} \in S$ with $\text{supp} \, f = \{m_1, \dots, m_n\}$ such that $\frac{c_j}d$ satisfies ACCP in $R$ for some $j \in \ldb 1,n \rdb$. After setting $g := \sum_{i=1}^n r_i X^{m_i}$, where $r_i = \frac{c_i}d$ for every $i \in \ldb 1, n \rdb$, it suffices to verify that~$g$ satisfies ACCP in $R[M]$. 
	
	Suppose, by way of contradiction, that there exists an ascending chain of principal ideals $(R g_n)_{n \in \nn_0}$ of $R[M]$ starting at $g_0 = g$ that does not stabilize. As mentioned in the introduction, since $M$ is torsion-free, we can assume that it has a total order compatible with the monoid operation and consider the degree of any nonzero element of $R[M]$ with respect to that order. Then we see that $(\deg g_n + M)_{n \in \nn_0}$ is an ascending chain of principal ideals of $M$, which must eventually stabilize because $M$ satisfies ACCP. Take $m \in \nn$ such that $\deg g_n = \deg g_m$ for every $n \ge m$. Now set $d_k : = \frac{g_{m+k-1}}{g_{m+k}}$ for every $k \in \nn$. Observe that $d_k \in R$ for every $k \in \nn$. In addition, for each $k \in \nn$, the equality $d_1 \cdots d_k g_{m+k} = g_m$, along with the fact that $g_m \mid_R g$, guarantees that $d_1 \cdots d_k \mid_R r_i$ for every $i \in \ldb 1,n \rdb$. Set $s_k := d_1 \cdots d_k$ for every $k \in \nn$, and observe that $(R \frac{r_j}{s_k})_{k \in \nn}$ is an ascending chain of principal ideals of $R$. This chain of principal ideals does not stabilize because $(R g_n)_{n \in \nn_0}$ does not stabilize. After setting $s_0 := 1$, we see that $(R \frac{r_j}{s_k})_{k \in \nn_0}$ is an ascending chain of principal ideals of $R$ starting at $r_j$ that does not stabilize, which contradicts that $r_j$ satisfies ACCP in $R$.
\end{proof}

In particular, the weak-ACCP property transfers from an integral domain to its polynomial ring.

\begin{cor} \label{cor:weak-ACCP and polynomials}
	For an integral domain $R$ the following statements are equivalent.
	\begin{enumerate}
		\item[(a)] $R$ is weak-ACCP.
		\smallskip
		
		\item[(b)] $R[X]$ is weak-ACCP.
		\smallskip
		
		\item[(c)] $R[X_\gamma \mid \gamma \in \Gamma]$ is weak-ACCP, where $\Gamma$ is a nonempty set.
	\end{enumerate}
\end{cor}

With notation as in Proposition~\ref{prop:weak-ACCP from monoid to monoid rings}, the constrain on the monoid $M$ to be reduced is not superfluous. This is illustrated in the following example.

\begin{example}
	Let $K$ be an algebraically closed field, and consider the monoid algebra $K[\qq]$. Notice that $\qq$ is torsion-free, and the fact that it is a group immediately implies that it satisfies ACCP. In addition, the fact that $K$ is a field ensures that it is a weak-ACCP domain. Arguing that $K[\qq]$ is not a weak-ACCP domain amounts to showing that it is antimatter. To do so, pick a nonunit $f \in K[\qq]^*$. After replacing~$f$ by one of its associates if necessary, we can assume that $f \in K[\qq_{\ge 0}]$ and $\text{ord} \, f = 0$. Now we can fix $n \in \nn$ large enough so that $f = g(X^{1/n})$ for some $g \in K[X]$ with $\deg g \ge 2$. As~$K$ is algebraically closed and~$g$ has degree at least $2$, one can write $g = g_1 g_2$ for some polynomials $g_1, g_2 \in K[X] \setminus K$. Observe that $\text{ord} \, f = 0$ implies that $\text{ord} \, g = 0$, and so $\text{ord} \, g_1 = \text{ord} \, g_2 = 0$. Thus, neither $g_1$ nor $g_2$ are monomials in $K[X]$, which implies that none of $g_1(X^{1/n})$ and $g_2(X^{1/n})$ is a unit in $K[\qq]$. Finally, the equality $f(X) = g_1(X^{1/n}) g_2(X^{1/n})$ ensures that $f$ is not irreducible in $K[\qq]$.
\end{example}

\begin{remark} \label{rem:strongly atomic domain that is not weak-ACCP}
	In \cite[Example~5.2]{mR93}, Roitman constructed an example of a strongly atomic domain $A$ such that $A[X]$ is not even atomic (the case when $k=2$). Then, in light of Corollary~\ref{cor:weak-ACCP and polynomials}, the integral domain $A$ is a strongly atomic domain that is not a weak-ACCP domain.
\end{remark}

With notation as in Proposition~\ref{prop:weak-ACCP from monoid to monoid rings}, one would like to replace the ACCP property on~$M$ by the weak-ACCP property. However, this is not possible, as the next example indicates.

\begin{example}
	Let $F$ be a field whose characteristic is not $2$, and let $(p_n)_{n \in \nn}$ be a strictly increasing sequence of primes with $p_1 > 3$. Then consider the additive submonoid of $\qq$
	\[
		M := \Big\langle \frac{1}{3^n p_n} \ \big{|} \ n \in \nn \Big\rangle.
	\]
	We have seen in Example~\ref{ex:Grams monoid is weak-ACCP} that $M$ is a weak-ACCP monoid. We proceed to show that the monoid algebra $F[M]$ is not weak-ACCP. To do so, consider the subset $S := \{X-1, X+1\}$ of $F[M]$, and let $f \in F[M]$ be a common divisor of $S$. Observe that $f(X) \mid_{F[M]} (X+1) - (X-1) = 2$, and so the fact that $M$ is a positive monoid ensures that $f \in F[M]^\times$. Observe that the sequence $\big( F[M](X^{1/3^n} - 1) \big)_{n \in \nn_0}$ is an ascending chain of principal ideals of $F[M]$ that starts at $X-1$ but does not stabilize; this is because
	\[
		X^{ \frac{1}{3^n} } - 1 = \big( X^{ \frac{1}{3^{n+1}} } - 1\big) \big( X^{ \frac{2}{3^{n+1}} } + X^{ \frac{1}{3^{n+1}} } + 1 \big).
	\]
	Similarly, we can verify that $\big( F[M]( X^{1/3^n} + 1) \big)_{n \in \nn_0}$ is an ascending chain of principal ideals of $F[M]$ starting at $X+1$ that does not stabilize. Hence neither $X-1$ nor $X+1$ satisfy ACCP in $F[M]$. Thus, $F[M]$ is not a weak-ACCP domain.
\end{example}

\smallskip
\subsection{Generalization of Grams' Construction}

One of our primary motivations to introduce the weak-ACCP property is that this notion allows us to construct atomic integral domains that do not satisfy ACCP.

\begin{theorem} \label{thm:Grams' construction and the weak-ACCP}
	Let $R$ be an integral domain, and let $M$ be a positive monoid. If $M$ is weak-ACCP, then so is the localization of $R[M]$ at the multiplicative subset $\{f \in R[M] \mid f(0) \neq 0\}$.
\end{theorem}

\begin{proof}
	Assume that $M$ is weak-ACCP. Then set $S := \{f \in R[M] \mid f(0) \neq 0\}$, and let $R[M]_S$ denote the localization of $R[M]$ at $S$. By virtue of Lemma~\ref{lem:associates weak-ACCP sets}, it suffices to argue that for any arbitrarily-chosen nonempty finite subset $T$ of $R[M]^*$, there is a common divisor $d'$ of $T$ in $R[M]_S$ that is a product of irreducibles such that~$\frac t{d'}$ satisfies ACCP in $R[M]_S$ for some $t \in T$. We can assume, without loss of generality, that $T$ and $R[M]_S^\times$ are disjoint. Set $E := \bigcup_{f \in T} \supp_{R[M]} \, f$. Since~$E$ is a nonempty finite subset of the weak-ACCP monoid~$M$, there is a common divisor $d$ of $E$ in $M$ such that $e - d$ satisfies ACCP in $M$ for some $e \in E$. Take $g \in T$ such that $e \in \text{supp}_{R[M]} \, g$, and set $f(X) := g(X)/X^d$. Because $d$ is a common divisor of $E$ in~$M$, it follows that $f \in R[M]$. In addition, after writing $f = \sum_{i=1}^k r_i X^{m_i}$ for some $r_1, \dots, r_k \in R$ with $r_1 \cdots r_k \neq 0$ and $m_1, \dots, m_k \in M$ with $m_1 > \dots > m_k$, we see that there is an index $j \in \ldb 1,k \rdb$ such that~$m_j$ satisfies ACCP in~$M$.
	
	Let us prove that $f$ satisfies ACCP in $R[M]_S$. Suppose, by way of contradiction, that there is an ascending chain of principal ideals $(R[M]_S f_n)_{n \in \nn_0}$ of $R[M]_S$ with $f_0 = f$ that does not stabilize. We can assume, without loss of generality, that $R[M]_S f_n \subsetneq R[M]_S f_{n+1}$ for every $n \in \nn_0$. Since $S \subseteq R[M]_S^\times$, we can further assume that $f_n \in R[M]$ for every $n \in \nn$. Now set
	\[
		b_0 := \min \big\{ m \in \supp_{R[M]} \, f \ \big{|} \ \, m \mid_M m_j \big\}.
	\]
	Since $m_j$ satisfies ACCP in $M$, the element $b_0$ also satisfies ACCP in $M$. Note that $m_j  \ge b_0 > 0$, where the second inequality comes from the fact that $f \notin R[M]_S^\times$. Take $g_1 \in R[M]$ and $s_1 \in S$ such that $s_1 f = f_1 g_1$. We claim that $b_0 \in \supp_{R[M]} \, s_1 f$. To argue this, take $m \in \supp_{R[M]} \, s_1$ and $i \in \ldb 1,k \rdb$ with $m + m_i = b_0$. Since $m_i \mid_M b_0$ and $b_0 \mid_M m_j$, it follows that $m_i \mid_M m_j$. This, along with the minimality of $b_0$, ensures that $m_i \ge b_0$. Thus, $m + m_i = b_0$ implies that $m_i = b_0$ (as~$M$ is a positive monoid), and so $m = 0$. As a result, the only monomial of degree $b_0$ we will obtain after multiplying out $s_1 f$ (and before cancellations) is $s_1(0) r_k X^{m_k}$, where $m_k = b_0$. Therefore we can take $c_1 \in \supp_{R[M]} \, f_1$ and $d_1 \in \supp_{R[M]} \, g_1$ such that $b_0 = c_1 + d_1$. Since none of $f_1$ and $g_1$ is a unit in $R[M]_S$, we see that $b_0 + M \subsetneq c_1 + M$. Now we can set
	\[
		b_1 := \min \big\{ m \in \supp_{R[M]} \, f_1 \, \big{|} \ \,  m \mid_M c_1 \big\},
	\]
	and note that $b_0 + M \subsetneq b_1 + M$. We can repeat this process indefinitely to obtain an ascending chain of principal ideals $(b_n + M)_{n \in \nn_0}$ of $M$ that does not stabilize, contradicting the fact that $b_0$ satisfies ACCP in~$M$.
\end{proof}

As illustrated in the following example, one can apply Theorem~\ref{thm:Grams' construction and the weak-ACCP} to produce integral domains that are weak-ACCP but do not satisfy ACCP.

\begin{example} \label{ex:weak-ACCP not ACCP}
	Fix a prime $p$, and let $(p_n)_{n \in \nn}$ be a strictly increasing sequence of primes that are different from $p$. In addition, let $q$ be a rational in $(0,1)$ such that $q^{-1} \notin \nn$. We have seen in Examples~\ref{ex:Grams monoid is weak-ACCP} and~\ref{ex:rational cyclic semirings is weak-ACCP but not ACCP} that the positive monoids
	\[
		G_p = \Big\langle \frac{1}{p^n p_n} \ \Big{|} \ n \in \nn_0 \Big\rangle \quad \text{ and } \quad S_q = \big\langle q^n \mid n \in \nn_0 \big\rangle
	\]
	are weak-ACCP but do not satisfy ACCP. Now fix a field $F$. It follows from Theorem~\ref{thm:Grams' construction and the weak-ACCP} that the integral domain $F[G_p]_S$, where $S$ is the multiplicative subset $\{f \in F[G_p] \mid f(0) \neq 0\}$ of $F[G_p]$, is a weak-ACCP domain, which does not satisfy ACCP because the ascending chain of principal ideals $\big( F[G_p]_S X^{1/p^n} \big)_{n \in \nn}$ does not stabilize. Similarly, we can deduce that $F[S_q]_T$, where $T$ is the multiplicative subset $\{f \in F[S_q] \mid f(0) \neq 0\}$ of $F[S_q]$ is a weak-ACCP domain that does not satisfy ACCP.
\end{example}
\medskip

At the end of Section~\ref{sec:algebraic constructions}, we have studied the behavior of the weak-ACCP property under directed unions. As another application of Theorem~\ref{thm:Grams' construction and the weak-ACCP}, we will show that the intersection of atomic (resp., weak-ACCP) domains may not be atomic (resp., weak-ACCP).

It is not hard to verify that the intersection of finitely many integral domains satisfying ACCP also satisfies ACCP. We can extend this result to certain type of non-necessarily finite intersections. For a nonempty set $\Gamma$, let $\{R_\gamma \mid \gamma \in \Gamma \}$ be a family of integral domains embedded into a common field. The integral domain $R = \bigcap_{\gamma \in \Gamma} R_\gamma$ is called the \emph{locally finite intersection} provided that for each $r \in R^*$, the set $\{\gamma \in \Gamma \mid r \notin R_\gamma^\times \}$ is finite. It is well known and not hard to argue that the locally finite intersection of integral domains satisfying ACCP also satisfies ACCP \cite[Corollary~4]{mZ88}. We proceed to argue that neither atomicity nor the weak-ACCP property are preserved under finite intersections.

Fix a field $F$, and let $(p_n)_{n \in \nn}$ and $(q_n)_{n \in \nn}$ be two strictly increasing sequences of odd primes such that $p_n < q_n < p_{n+1}$ for every $n \in \nn$. Now consider the positive monoids
\[
	G_1 := \Big\langle \frac{1}{2^n p_n} \ \Big{|} \ n \in \nn \Big\rangle \quad \text{and} \quad G_2 := \Big\langle \frac{1}{2^n q_n} \ \Big{|} \ n \in \nn \Big\rangle.
\]
Let $L_1$ and $L_2$ be the respective localizations of $F[G_1]$ and $F[G_2]$ at their multiplicative subsets consisting of polynomial expressions of order $0$. Observe that we have constructed $L_1$ and $L_2$ mimicking Grams' construction, except that we have slightly changed the exponent monoids. Therefore $L_1$ and~$L_2$ are slightly modified versions of the Grams monoid, and so it follows from Theorem~\ref{thm:Grams' construction and the weak-ACCP} that they are both weak-ACCP. The integral domains $L_1$ and $L_2$ are both inside $F(\qq)$, the field of fractions of $F[\qq]$.
With the notation introduced in this paragraph, we proceed to prove the following result.

\begin{prop} \label{prop:antimatter intersection of two weak-ACCP domains}
	The integral domain $L_1 \cap L_2$ is antimatter.
\end{prop}

\begin{proof}
	Set $R := L_1 \cap L_2$. Observe that $R^\times = L_1^\times \cap L_2^\times$. An elementary divisibility argument can be used to argue that $M := G_1 \cap G_2 = \big\langle \frac{1}{2^n} \mid n \in \nn \big\rangle$. In addition, an argument completely similar to that developed in \cite[Lemma~1.1]{aG74} shows that $M$ is a greatest-divisor submonoid of both $G_1$ and $G_2$ (see Example~\ref{ex:Grams monoid has a greatest-divisor submonoid}). Then for every $f \in F[G_1] \cup F[G_2]$, we can define
	\[
	\mu(f) := 
		\begin{cases}
			\min \{ \text{gd}_M(\alpha) \text{ in } G_1 \mid \alpha \in \text{supp} \, f \} &  \text{ if } \ f \in F[G_1] \\
			\min \{ \text{gd}_M(\alpha) \text{ in } G_2 \mid \alpha \in \text{supp} \, f \} &  \text{ if } \ f \in F[G_2].
		\end{cases}
	\]
	Observe that $\mu$ is well defined in $F[M] = F[G_1] \cap F[G_2]$ as for each $\alpha \in M$ the equality $\gd_M(\alpha) = \alpha$ holds in both $G_1$ and $G_2$. We proceed to prove the following claim.
	\smallskip
	
	\noindent \textit{Claim.} $R = \big\{ X^q u(X) \mid q \in M \text{ and } u(X) \in L_1^\times \cap L_2^\times \big\}$.
	\smallskip
	
	\noindent \textit{Proof of Claim.} Observe that 
	$\big\{ X^q u(X) \mid q \in M \text{ and } u(X) \in L_1^\times \cap L_2^\times \big\}$ is contained in $R$. To argue the reverse implication, fix an element in $R$, namely, $\frac{f_1}{s_1} = \frac{f_2}{s_2} \in R$, where $f_1, s_1 \in F[G_1]$ with $s_1(0) \neq 0$ and $f_2, s_2 \in F[G_2]$ with $s_2(0) \neq 0$. Let us argue that $\mu(f_1) = \mu(f_2)$. Suppose, by way of contradiction, that $\mu(f_1) \neq \mu(f_2)$ and assume, without loss of generality, that $\mu(f_1) < \mu(f_2)$. This implies that
	\[
		s_2(X) \frac{f_1(X)}{X^{\mu(f_1)}} = X^{\mu(f_2) - \mu(f_1)} s_1(X) \frac{f_2(X)}{X^{\mu(f_2)}}.
	\]
	Now set $r := \min \big\{ \alpha \in \text{supp} \, f_1(X)/X^{\mu(f_1)} \mid \text{gd}_M(\alpha) = 0 \big\}$. Some elementary arguments can be used to verify that $r \in \text{supp} \, s_2(X) f_1(X)/X^{\mu(f_1)}$. On the other hand, some elementary arguments can also be used to verify that $\mu(f_2) - \mu(f_1) \nmid_{G_1 + G_2} r$, which implies that $r \notin \text{supp} \, X^{ \mu(f_2) - \mu(f_1)} s_1(X) f_2(X)/X^{\mu(f_2)}$, which is a contradiction. Hence $\mu(f_1) = \mu(f_2)$, as desired.
	
	Now set $q = \mu(f_1)$ and then set $r := \min \{ \alpha \in \text{supp} \, f_1(X)/X^q \mid \text{gd}_M(\alpha) = 0 \}$. We split the rest of our argument into the following two cases.
	\smallskip
	
	\textsc{Case 1:}  $r=0$. Write $f_1(X)/X^q := c + g(X)$ for some $c \in F$ and $g \in F[G_1]$ with $\text{ord} \, g > 0$. The fact that $r=0$ implies that $c \neq 0$ and, therefore, we see that $(f_1(X)/X^q)/s_1(X) \in L_1^\times$. In addition, from the identity $s_2(X) (c + g(X)) = s_1(X) f_2(X)/X^q$, we deduce that $f_2(X)/X^q$ has constant coefficient $c s_2(0)/s_1(0)$, which is different from $0$. Thus, $(f_1(X)/X^q)/s_1(X)  = (f_2(X)/X^q)/s_2(X) \in L_2^\times$. As a consequence, $\frac{f_1}{s_1} \in \{ X^q u(X) \mid q \in M \text{ and } u(X) \in L_1 \cap L_2 \}$, and so the reverse inclusion also holds in this case.
	\smallskip
	
	\textsc{Case 2:} $r > 0$. Set $h(X) := f_2(X)/X^q$, and observe that $r > 0$ implies that $\text{ord} \, h > 0$. It is not hard to verify that $r \in \text{supp} \, s_2(X)(f_1(X)/X^q)$. Observe that $\alpha \nmid_{G_1 + G_2} r$ for any $\alpha \in M^\bullet$. In addition, one can check that $\frac{1}{2^k q_k} \nmid_{G_1 + G_2} r$ for any $k \in \nn$. Combining the last two statements, we obtain that $r \notin \text{supp} \, s_1(X)(f_2(X)/X^q)$, which is a contradiction. Hence the claim is established.
	\smallskip
	
	Finally, let $h$ be a nonzero nonunit in $R$. The established claim allows us to write $h(X) = X^q u(X)$ for some $q \in M^\bullet$ and $u \in R^\times = L_1^\times \cap L_2^\times$.	Observe that $X^q \notin R^\times$ because $X^q \notin L_1^\times$. Similarly, $X^{q/2}$ and $X^{q/2}u(X)$ are nonunits in $R$. Therefore the decomposition $h(X) = X^{q/2}(X^{q/2}u(X))$ shows that~$h$ is not irreducible. Hence $R$ is antimatter.
\end{proof}

Proposition~\ref{prop:antimatter intersection of two weak-ACCP domains} allows us to conclude with the following remark.

\begin{remark}
	The intersection of two weak-ACCP (resp., atomic) domains is not necessarily a weak-ACCP (resp., atomic) domain.
\end{remark}

\smallskip
\subsection{A Finite-Rank Atomic Monoid Algebra without ACCP}

As mentioned in the introduction, every known atomic monoid algebra that does not satisfy ACCP has infinite rank. In this subsection, we construct a rank-$2$ monoid algebra satisfying the weak-ACCP property but not the ACCP property.
\smallskip

Let $(p_n)_{n \in \nn}$ be a strictly increasing sequence of odd primes such that $\sum_{n=1}^\infty \frac 1{p_n} < \frac{1}{3}$, and consider the positive monoid $P = \big\langle \frac{1}{p_n} \mid n \in \nn \big\rangle$. It is not hard to check that every element $q \in P$ can be written in an essentially unique way as
\begin{equation} \label{eq:unique decomposition}
	q = n_0 + \sum_{i=1}^\ell \frac{n_i}{p_i}
\end{equation}
where $n_0, n_1, \dots, n_\ell \in \nn_0$ and $n_i \in \ldb 0, p_i - 1 \rdb$ for every $i \in \ldb 1, \ell \rdb$. From the uniqueness of such a representation, we can deduce that $P$ satisfies ACCP (see \cite[Example~2.1]{AAZ90} and \cite[Proposition~4.2]{fG22}). Now consider the following subset of $P$:
\[
	A := \bigg\{ \frac 1{p_{j_k}} + \sum_{i=1}^\ell \frac{1}{p_{j_i}}  \ \Big{|} \ k, \ell \in \nn \text{ with } k \in \ldb 1, \ell \rdb \text{ and } j_1 < j_2 < \cdots < j_\ell \bigg\}.
\]
In addition, let $\beta$ be an irrational number such that $\beta > 1$, and consider the following subset of $\nn_0 \beta + \qq$:
\[
	B := \{ \beta \} \bigcup \bigg\{ \beta_\ell :=  \beta - \sum_{i=1}^\ell \frac{1}{p_i}  \ \Big{|} \ \ell \in \nn \bigg\}.
\]
For convenience of notation, set $\beta_0 = \beta$. For the rest of this section, fix a field $F$ and let $M$ denote the positive monoid generated by $A \cup B$. We will prove that the monoid algebra $F[M]$ is a weak-ACCP domain that does not satisfy ACCP.

\begin{prop} \label{prop:atomic monoid that is not weak-ACCP}
	The monoid $M$ is atomic but does not satisfy ACCP.
\end{prop}

\begin{proof}
	It suffices to show that $\mathcal{A}(M) = A \cup B$. Since $\beta$ is irrational, none of the elements in $A$ is divisible in $M$ by any element of $B$. Therefore $\langle A \rangle$ is a divisor-closed submonoid of~$M$. Now take $a \in A$ and write $a = a_1 + \dots + a_k$ for some $a_1, a_2, \dots, a_k \in A$. As in each sum defining the elements of $A$ only one atom of $P$ repeats and repeats exactly twice, $k=1$ (here we used that each $p_n$ is odd). This, along with the fact that $\langle A \rangle$ is a divisor-closed submonoid of $M$, ensures that $a \in \mathcal{A}(M)$. Now fix $b \in B$. Since the set $\{1, \beta\}$ is linearly independent over $\qq$, none of the elements in $2B + M$ can divide $b$ in $M$. Write $b = b' + a'_1 + \dots + a'_\ell$ for some $b' \in B$ and $a'_1, \dots, a'_\ell \in A$. Observe that $b - b' \in P$ and also that the decomposition~\eqref{eq:unique decomposition} of $b-b'$ does not repeat any atom of $P$. This guarantees that $\ell = 0$, and so $b \in \mathcal{A}(M)$. Hence $\mathcal{A}(M) = A \cup B$. Since the sequence $(2 \beta_n + M)_{n \in \nn}$ is an ascending chain of principal ideals of $M$ that does not stabilize, we conclude that $M$ is an atomic monoid that does not satisfy ACCP. 
\end{proof}

\begin{lemma}  \label{lem:sufficient condition for belongin to <A>}
	For $n_1, \dots, n_\ell \in \nn_0$, set $q = \sum_{i=1}^\ell \frac{n_i}{p_i}$. If $\min\{n_j, n_k\} \ge 2$ for some $j,k \in \ldb 1, \ell \rdb$, then $q \in \langle A \rangle$.
\end{lemma}

\begin{proof}
	The argument follows, \emph{mutatis mutandis}, that given in the proof of \cite[Lemma~4.2]{GL23}.
\end{proof}

Now define the function $\varphi \colon M \to \nn_0$ by letting $\varphi(r)$ be the unique $n_0 \in \nn_0$ such that $r = n_0 \beta + q$ for some $q \in \qq$.

\begin{lemma} \label{lem:auxiliary for main theorem}
	If $\varphi(r) \ge 2$ for some $r \in M$, then there exists $N \in \nn$ such that $\varphi(r) \beta_n \mid_M r$ for every $n \in \nn$ with $n \ge N$.
\end{lemma}

\begin{proof}
	Observe that we can write $r = \varphi(r) \beta + \sum_{i=1}^\ell \frac{n_i}{p_i}$ for some $n_1, \dots, n_\ell \in \zz_{\ge -\varphi(r)}$. Fix now $N \in \nn$ with $N \ge \ell + 2$. Since
	\[
		r - \varphi(r) \beta_N =  \sum_{i=1}^\ell \frac{\varphi(r) + n_i}{p_i} + \sum_{i = \ell + 1}^N \frac{\varphi(r)}{p_i},
	\]
	it follows from Lemma~\ref{lem:sufficient condition for belongin to <A>} that $\varphi(r) \beta_N \mid_M r$. Hence $\varphi(r) \beta_n \mid_M r$ for every $n \in \nn$ with $n \ge N$.
\end{proof}

\begin{lemma} \label{lem:representation of elements in B + <A>}
	If $r \in B + \langle A \rangle$, then $r$ can be written as $r = \beta + n + \sum_{i=1}^\ell \frac{n_i}{p_i}$ for some $\ell, n \in \nn_0$ and $n_1, \dots, n_\ell \in \zz$ with $n_i \in \ldb -1, p_i - 2 \rdb$ for every $i \in \ldb 1, \ell \rdb$. Moreover, in such a decomposition, $n$ and the coefficients $n_1, \dots, n_\ell$ are uniquely determined.
\end{lemma}

\begin{proof}
	Take $r \in B + \langle A \rangle$. Note that we can write $r = \beta + n + \sum_{i=1}^\ell \frac{n_i}{p_i}$ for some $\ell, n \in \nn_0$ and $n_1, \dots, n_\ell \in \zz$ with $n_i \in \ldb -1, p_i - 2 \rdb$ for every $i \in \ldb 1, \ell \rdb$. The uniqueness of such a decomposition can be argued similarly to that of~\eqref{eq:unique decomposition}.
\end{proof}

Now define $\psi \colon B + \langle A \rangle \to P$ as follows. For $r \in B + \langle A \rangle$, let $\psi(r) = n + \sum_{i=1}^\ell \frac{\max\{n_i,0\}}{p_i}$, where $r = \beta + n + \sum_{i=1}^\ell \frac{n_i}{p_i}$ is the unique decomposition described in the statement of Lemma~\ref{lem:representation of elements in B + <A>}.

\begin{lemma} \label{lem:function psi}
	If $a \in A$ and $b \in B + \langle A \rangle$, then the following statements hold.
	\begin{enumerate}
		\item $\psi(b) \neq \psi(b+a)$.
		\smallskip
		
		\item $\psi(b) \mid_P \psi(b+a)$.
	\end{enumerate}
\end{lemma}

\begin{proof}
	Write $b = \beta + n + \sum_{i=1}^\ell \frac{n_i}{p_i}$ for some $\ell, n \in \nn_0$ and $n_1, \dots, n_\ell \in \zz$ with $n_i \in \ldb -1, p_i - 2 \rdb$ for every $i \in \ldb 1, \ell \rdb$. By virtue of the uniqueness of the decomposition~\eqref{eq:unique decomposition}, there exists a unique $k \in \nn$ such that $\frac{2}{p_k} \mid_P a$. By inserting zero coefficients $n_j$ if necessary in the decomposition of $b$, we can assume that $k \le \ell$. Now write $b + a= \beta + n' + \sum_{i=1}^\ell \frac{n'_i}{p_i}$ for some $\ell, n' \in \nn_0$ and $n'_1, \dots, n'_\ell \in \zz$ with $n'_i \in \ldb -1, p_i - 2 \rdb$ for every $i \in \ldb 1, \ell \rdb$
	\smallskip
	
	(1) If $n_k = p_k - 3$ (resp., $n_k = p_k - 2$), then $n'_k = -1$ (resp., $n'_k = 0$) and so $\max\{n_k, 0\} \neq \max \{n'_k, 0\}$ (observe that $3 \notin \{p_n \mid n \in \nn\}$ because $\sum_{n=1}^\infty \frac 1{p_n} < \frac{1}{3}$). On the other hand, if $n_k \in \ldb -1, p_k - 4 \rdb$, then $\max \{n_k, 0\} \le  n_k + 1 < n'_k = \max\{n'_k,0\}$. We conclude that $\max\{n_k, 0\} \neq \max\{n'_k, 0\}$ and, therefore, $\psi(b) \neq \psi(b+a)$.
	\smallskip
	
	(2) Observe that $\psi$ can be naturally extended to a map $B + P \to P$, which we also call $\psi$. Thus, it suffices to argue that $\psi(b) \mid_P \psi\big(b + \frac 1{p_j}\big)$ for every $j \in \nn$. Fix $j \in \nn$. As we did for $k$, we can assume that $j \le \ell$. If $n_j = -1$, then $\psi(b + \frac{1}{p_j}) = \psi(b)$. In addition, if $n_j \in \ldb 0, p_j - 3 \rdb$, then $\psi(b + \frac 1{p_j}) = \psi(b) + \frac 1{p_j}$. Finally, if $n_j = p_j - 2$, then $\psi(b + \frac 1{p_j}) = b + \frac 2{p_j}$. In all possible cases, we see that $\psi(b) \mid_P \psi(b + \frac 1{p_j})$, which concludes the proof.
\end{proof}

We are in a position to prove our last main result.

\begin{theorem} \label{thm:weak-ACCP monoid ring not ACCP}
	The monoid algebra $F[M]$ is weak-ACCP but does not satisfy ACCP.
\end{theorem}

\begin{proof}
	Set $R := F[X;\qq]$. Now we identify $F[Y;M] := F[M]$ with a subring of $R[x]$ by using the natural embedding of $F[M]$ into $R[x]$ given by the assignments $Y^{n_0 \beta + q} \mapsto X^q x^{n_0}$ for any $n_0 \in \nn_0$ and $q \in \qq$. For $f \in R[x]$, we let $\deg f$ (resp., $\text{ord} \, f$) stand for the degree (resp., the order) of $f$ with respect to the indeterminate~$x$. If $f \in R$, then the degree (resp., order) of $f$ with respect to the variable $X$ is denoted by $\deg_R f$ (resp., $\text{ord}_R \, f$). Observe that $R \cap F[M] \subseteq F[P]$. 
	
	To prove that $F[M]$ is a weak-ACCP domain, let $S$ be a nonempty finite subset of $F[M]^*$. Since $\inf \langle B \rangle^\bullet \ge \frac 12$, we can take $b \in \langle B \rangle$ such that $Y^b \mid_{F[M]} s$ for all $s \in S$ and such that no common divisor of the set $\{sY^{-b} \mid s \in S\}$ has the form $Y^{\beta_n}$ for any $n \in \nn_0$. Hence Lemma~\ref{lem:auxiliary for main theorem} ensures the existence of $f \in \{sY^{-b} \mid s \in S\}$ with $\text{ord} \, f < 2$. Because $M$ is atomic by Proposition~\ref{prop:atomic monoid that is not weak-ACCP}, we see that $Y^b$ factors into irreducible in $F[M]$, and so we are done once we prove the following claim.
	\smallskip
	
	\noindent \textit{Claim.} $f$ satisfies ACCP.
	\smallskip
	
	\noindent \textit{Proof of Claim.} Suppose, by way of contradiction, that there is an ascending chain of principal ideals $(F[M] f_n)_{n \in \nn_0}$ of $F[M]$ that does not stabilize with $f_0 = f$. We can further assume that $F[M] f_n \subsetneq F[M] f_{n+1}$ for every $n \in \nn_0$. Take a sequence $(g_n)_{n \in \nn}$ whose terms belong to $F[M]$ such that $f_n = f_{n+1} g_{n+1}$ for every $n \in \nn_0$. Since the sequence $(\deg f_n)_{n \in \nn}$ of nonnegative integers eventually stabilizes, we can fix $K \in \zz_{\ge 2}$ 
	such that $g_n \in R \cap F[M] \subseteq F[P]$ for every $n \ge K$. For each $n \ge K$, the fact that $F[M] f_n \subsetneq F[M] f_{n+1}$ implies that $g_n \in F[P] \setminus F$. We split the rest of the proof into the following two cases.
	\smallskip
	
	\textsc{Case 1:} $\text{ord} \, f = 0$. Since the positive monoid $P$ satisfies ACCP, the monoid algebra $F[P]$ also satisfies ACCP (see \cite[Proposition~4.2]{fG22}). For each $n \ge K$, the fact that $g_n \in F[P]$ ensures that $\deg \, g_n = 0$, and so $f_n(0)/f_{n+1}(0) = g_{n+1} \in (F[P] \setminus F) \cap F[M]$. Therefore $(F[P] f_n(0))_{n \ge K}$ is an ascending chain of principal ideals in $F[P]$. Now the fact that $(F[P] f_n(0))_{n \ge K}$ does not stabilize contradicts that $F[P]$ satisfies ACCP.
	\smallskip
	
	\textsc{Case 2:} $\text{ord} \, f = 1$. If $\text{ord} \, f_m = 0$ for some $m \in \nn$, then we can proceed as in Case~1 to generate a contradiction. Therefore we assume that, for each $n \in \nn_0$, the equality $\text{ord} \, f_n = 1$ holds (that is, $f'_n(0) \neq 0$, where $f'_n$ is the formal derivative of $f'_n$ in $R[x]$). For every $n \in \nn$, the equality $f_K = f_{K+n} \prod_{i=1}^n g_{K+i}$ holds, and so $f'_K(0) = f'_{K+n}(0) \prod_{i=1}^n g_{K+i}$. Now set
	\[
		d_n = \deg_R g_n \in \langle A \rangle^\bullet \quad \text{ and } \quad b_n := \beta + \deg_R f'_n(0) \in B + \langle A \rangle
	\]
	for every $n \ge K$. Therefore $b_k = b_{K+n} + \sum_{i=1}^n d_i$ for every $n \in \nn$ and, by virtue of Lemma~\ref{lem:function psi}, the sequence $\big( \psi(b_n) + P \big)_{n \ge K}$ is an ascending chain of principal ideals of $P$ that does not stabilize. However, this contradicts that $P$ satisfies ACCP.
\end{proof}

\smallskip
\subsection{An Atomic Monoid Algebra without ACCP by Zaks}
 
 As far as we know, the monoid algebra we have constructed in the previous subsection is the first known \emph{finite-rank} atomic monoid algebra that does not satisfy ACCP. However, infinite-rank atomic monoid algebras that do not satisfy ACCP have been constructed by Zaks in~\cite{aZ82} and more recently by the authors in~\cite{GL23}. We conclude this section with a few words about the construction given by Zaks. Although it is not clear \emph{a priori} that the integral domain in Zaks' construction is a monoid algebra, this can be seen after some algebraic manipulations, and we do so in this subsection. We proceed to remind the description of the integral domain constructed by Zaks. 

\begin{example} \label{ex:Zaks domain}
	Let $F$ be a field, and consider the set of indeterminates $U, V, W$, and $X_n$ for all $n \in \nn$ over $F$. Now set $ Y_n := UV / (W^n X_n)$ for every $n \in \nn$, and let $R$ be the smallest subring of $F[U, V, W^{\pm 1}, X_n^{\pm 1}]_{n \in \nn}$ containing the set $\{U, V, W, X_n, Y_n \mid n \in \nn \}$. It was proved by Zaks in~\cite{aZ82} that~$R$ is an atomic domain that does not satisfy ACCP.
\end{example}

Motivated by Example~\ref{ex:Zaks domain}, we make the following definition.

\begin{definition}
	We call the integral domain $R$ in Example~\ref{ex:Zaks domain} the \emph{Zaks domain} over the field $F$, or a \emph{Zaks domain} for short.
\end{definition}

We proceed to verify that Zaks domains are monoid algebras.

\begin{prop}
	Every Zaks domain is an infinite-rank monoid algebra.
\end{prop}

\begin{proof}
	Let $M$ be the additive monoid $\nn \oplus \nn \oplus \zz \oplus \zz^{\oplus \nn}$. With notation as in Example~\ref{ex:Zaks domain}, we can think of the ring $F[U, V, W^{\pm 1}, X_n^{\pm 1}]$ as the monoid algebra $F[M;Z]$, where $Z^m = U^a V^b W^c (\prod_{n \in \nn} X_n^{d_n})$ for every element $m = (a,b,c, (d_n)_{n \in \nn}) \in M$ (here all but finitely many entries $d_n$ are zero). As a result, the Zaks domain $R$ over $F$ is the monoid algebra $F[N;Z]$, where $N$ is the submonoid of $M$ generated by the elements $(1,0,0,0, \ldots)$, $(0,1,0,0, \ldots)$, $(0,0,1,0, \ldots)$, $x_n := (0,0,0,\dots, 0,1,0, \dots)$, and $y_n := (-n, 1, 1, 0, 0, \ldots, 0, -1, 0, \ldots)$, where the entry $1$ of $x_n$ and the entry $-1$ of $y_n$ appear in the $(n+3)$-th position for every $n \in \nn$. Thus, $R$ is indeed a monoid algebra. Finally, it is clear that the monoid $N$ has infinite rank, so $R$ is an infinite-rank monoid algebra.
\end{proof}

We conclude this section with the following question.

\begin{question}
	Is every Zaks domain a weak-ACCP domain?
\end{question}

\bigskip
\section{A Weak-ACCP Pullback Ring}
\label{sec:pullback}

In this final section, we discuss an example of an atomic domain not satisfying ACCP that was given by Boynton and Coykendall in~\cite{BC19} as a pullback of certain commutative rings. We will prove that such an integral domain is weak-ACCP. Let $D$ be an integral domain with quotient field  $K$, and take distinct irreducible polynomials $a_1(X), \dots, a_n(X) \in K[X]$. For each $i \in \ldb 1, n \rdb$, fix a root $\theta_i$ of $a_i(X)$ and an overring $D_i$ of $D[\theta_i]$. Now set $f(X) = a_1(X) \cdots a_n(X)$ and $I := K[X] f(X)$. The canonically defined homomorphisms
\[
	\eta \colon K[X] \to K[X]/I \cong \prod_{i=1}^{n} K[\theta_i] \quad \text{ and } \quad \iota \colon \prod_{i=1}^n D_i \hookrightarrow \prod_{i=1}^n K[\theta_i]
\]
are surjective and injective, respectively. One can easily see that the pullback ring of these two homomorphisms can be written as
\[
	R := \big\{ g(X) \in K[X] \mid g(\theta_i) \in D_i \text{ for every } i \in \ldb 1,n \rdb \big\}.
\]
The atomicity of $R$ was characterized in~\cite[Theorem~3.3]{BC19}, and this result was used to give the following example of an atomic domain that does not satisfy ACCP.
\smallskip

\begin{example} \label{ex:BC pullback example}
	Let $F$ be a field, and let $x$ and $y$ be two distinct indeterminates. Consider the integral domain $D := S^{-1}F[x,y]$ that we obtain after localizing $F[x,y]$ at its multiplicative subset
	\[
		S := \{r \in F[x,y] \mid r \notin (x) \text{ and } r \notin (y)\}.
	\]
	With the choices given by the equalities $a_1(X) := X$ and $a_2(X) := X-1$ as well as $D_1 := D[\frac{y}{x^2}]$ and $D_2 := D[\frac{x}{y^2}]$ in the pullback construction defined above, we obtain a pullback domain $R$ that is atomic but does not satisfy ACCP (see \cite[Example~3.5]{BC19}).
\end{example}

The integral domain $R$ in Example~\ref{ex:BC pullback example} is the primary object of this section, and so it is convenient to give it a name.

\begin{definition}
	We call the integral domain $R$ in Example~\ref{ex:BC pullback example} the \emph{Boynton-Coykendall domain} over the field $F$, or a \emph{BC-domain} for short.
\end{definition}

 We proceed to refine the fact that the BC-domain is atomic, by proving that it is indeed a weak-ACCP domain. It is worth emphasizing that we will not use the fact already established by Boynton and Coykendall~\cite{BC19} that the BC-domain is atomic. In addition, the proof that we offer here is neither based on nor similar to that given in~\cite{BC19}. We need the following lemma.

\begin{lemma} \label{lem:pullback auxiliar}
	With notation as in Example~\ref{ex:BC pullback example}, let $f$ be an element of the BC-domain and let $K$ be the quotient field of~$D$. Then the following conditions are equivalent.
	\begin{enumerate}
		\item[(a)] $f$ does not satisfy ACCP.
		\smallskip
		
		\item[(b)] $f(0) = f(1) = 0$.
		\smallskip
		
		\item[(c)] $X(X-1) \mid_{K[X]} f(X)$.
	\end{enumerate}
\end{lemma}

\begin{proof}
	Note that 
	$D^\times = \big\{ \frac{s_1}{s_2} \mid s_1, s_2 \in S \big \}$. Observe that every nonzero element $p \in D \setminus D^\times$ can be written as $\frac rs$ for some $s \in S$ and $r \in F[x,y]$ such that $r$ is divisible by either $x$ or $y$ in $F[x,y]$. After factoring from $r$ as many copies of both $x$ and $y$ as possible, we see that the equality $p = ux^jy^k$ holds for some $u \in D^\times$ and $j,k \in \nn_0$ such that $j+k \ge 1$. In addition, note that $D_1^\times = D_2^\times = D^\times$. Finally, $R^\times = K \cap D_1^\times \cap D_2^\times = D^\times$.
	\smallskip
	
	(a) $\Rightarrow$ (b): Suppose, by way of contradiction, that $f(0) \neq 0$. Since $f$ does not satisfy ACCP, we can take a strictly ascending chain of principal ideals $(R f_n)_{n \in \nn_0}$ with $f_0 = f$. For each $n \in \nn_0$, write $f_n = g_{n+1} f_{n+1}$ for some $g_{n+1} \in R$. Since the sequence $(\deg f_n)_{n \in \nn}$ is decreasing, there exists $N \in \nn$ such that $\deg g_n = 0$ for every $n \ge N$. This, along with the fact that $D_1^\times = D^\times = R^\times$, guarantees that $(D_1 f_{N+k}(0))_{k \in \nn_0}$ is a strictly ascending chain of principal ideals in $D_1$. Because $D$ is obtained by localizing the Noetherian ring $F[x,y]$, we see that $D$ is also a Noetherian ring. In addition, $D_1 = D[\frac{y}{x^2}]$ is a homomorphic image of the polynomial ring $D[X]$, which is Noetherian by virtue of Hilbert Basis Theorem. Hence $D_1$ is a Noetherian ring. Therefore the chain $(D_1 f_{N+k}(0))_{k \in \nn_0}$ must stabilize, a contradiction. We can proceed similarly to produce a contradiction under the assumption that $f(1) \neq 0$. Hence $f(0) = f(1) = 0$.
	\smallskip
	
	(b) $\Rightarrow$ (c): This is straightforward.
	\smallskip
	
(c) $\Rightarrow$ (a): Finally, assume that $X(X-1) \mid_{K[X]} f(X)$. For each $n \in \nn$, the fact that $x^n$ divides $X(X-1)$ in $R$ now implies that $\big( R \frac{f}{x^n} \big)_{n \in \nn_0}$ is an ascending chain of principal ideals in $R$ starting at~$f$. Since $x \notin R^\times$, the same chain of principal ideals does not stabilize, and so $f$ does not satisfy ACCP.
\end{proof}

We are in a position to prove that the BC-domain is a weak-ACCP domain and, therefore, to identify another example of a weak-ACCP domain that does not satisfy ACCP (cf. Example~\ref{ex:weak-ACCP not ACCP} and Theorem~\ref{thm:weak-ACCP monoid ring not ACCP}).

\begin{theorem}
	Every BC-domain is a weak-ACCP domain.
\end{theorem}

\begin{proof}
	For this proof we adopt the notation introduced in Example~\ref{ex:BC pullback example}. In particular, $R$ denotes the BC-domain over the field $F$. Let $W$ be a nonempty finite subset of $R$. If there exists $w \in W$ such that $X(X-1)$ does not divide $w(X)$ in $K[X]$, then $1$ is a common divisor of $W$ and $\frac w1$ satisfies ACCP by Lemma~\ref{lem:pullback auxiliar}. Therefore we will assume for the rest of the proof that $X(X-1)$ divides each element of~$W$ in $K[X]$. Now consider
	\[
		m:= \max \big\{ n \in \nn_0 \mid X^n(X-1)^n \mid_{K[X]} w(X) \text{ for all } w \in W \big\}.
	\]
	For each $w \in W$, write $w(X) = c_w X^m(X-1)^m h_w(X)$ for some $c_w \in K$ and $h_w(X) \in D[X] \subset R$. Then there exists $N \in \nn$ such that $x^N y^N c_w \in D$ for every $w \in W$. We claim that
	\[
		d(X) := \left( \frac{x}{y^2} \right)^N X^m \bigg( \frac{y}{x^2} \bigg)^N (X-1)^m
	\]
	is a common divisor of $W$ that factors into irreducibles such that $\frac wd$ satisfies ACCP for some $w \in W$. Since none of the polynomials $( \frac{x}{y^2})^N X^m$ and $(\frac{y}{x^2})^N (X-1)^m$ contains the set $\{0,1\}$ inside its set of roots, it follows from Lemma~\ref{lem:pullback auxiliar} that both $( \frac{x}{y^2})^N X^m$ and $(\frac{y}{x^2})^N (X-1)^m$ satisfy ACCP. As a consequence, $d$ must factor into irreducibles. Now observe that
	\[
		\frac{w(X)}{d(X)} = x^Ny^N c_w h_w(X) \in D[X] \subseteq R
	\]
	for every $w \in W$, and so $d$ is a common divisor of $W$. Finally, it follows from the maximality of $m$ that there exists $w \in W$ such that $X(X-1)$ does not divide $w(X)/d(X)$ in $K[X]$. Hence it follows from Lemma~\ref{lem:pullback auxiliar} that $\frac wd$ satisfies ACCP in~$R$. We conclude that $R$ is a weak-ACCP domain.
\end{proof}

As a corollary, we obtain the conclusion of \cite[Example~3.5]{BC19}.

\begin{cor}
	Every BC-domain is atomic.
\end{cor}

\bigskip
\section*{Acknowledgments}

During the preparation of this paper, both authors were part of PRIMES-USA at MIT, and they would like to thank the directors and organizers of the program for making this collaboration possible. Finally, the first author kindly acknowledges support from the NSF under the awards DMS-1903069 and DMS-2213323.

\bigskip

\end{document}